\documentclass[aap,MSNbibl,seceqn,dvips]{arximspdf}

\doi{10.1214/11-AAP838} 
\volume{23}
\issue{1}
\pubyear{2013}
\firstpage{251}
\lastpage{282}

\makeatletter
\newcommand{\eqref}[1]{(\ref{#1})}

\newcommand{\dd}{\mathrm{d}}
\newcommand{\BE}{\mathbf{E}}
\newcommand{\BP}{\mathbf{P}}
\newcommand{\R}{\mathbb{R}}
\newcommand{\argmax}{\operatorname{argmax}}
\renewcommand{\geqslant}{\geq}
\renewcommand{\leqslant}{\leq}

\newtheorem{thmm}{Theorem}[section]
\newproclaim{defn}[section]{Definition}

\newtheorem{coro}[thmm]{Corollary}
\newtheorem{lemma}[thmm]{Lemma}
\newproclaim{assmp}{Assumption}
\newproclaim{exam}{Example}
\newproclaim{remark}{Remark}
\newproclaim{com}{Comments}

\makeatother

\begin{document}
\begin{frontmatter}

\title{Optimal stopping under probability distortion}
\runtitle{Optimal stopping under probability distortion}

\begin{aug}
\author[A]{\fnms{Zuo Quan} \snm{Xu}\ead[label=e1]{maxu@inet.polyu.edu.hk}\thanksref{t1,t3}}
\and
\author[B]{\fnms{Xun Yu} \snm{Zhou}\corref{}\ead[label=e2]{zhouxy@maths.ox.ac.uk}\thanksref{t2,t3}}
\thankstext{t1}{Supported by a start-up fund of the Hong Kong
Polytechnic University.}
\thankstext{t2}{Supported by a
start-up fund of the University of Oxford.}
\thankstext{t3}{Supported by grants from the Nomura Centre for Mathematical
Finance and the Oxford--Man Institute of Quantitative Finance.}
\runauthor{Z.~Q. Xu and X.~Y. Zhou}
\affiliation{Hong Kong Polytechnic University, and
University of Oxford and~Chinese~University of Hong Kong}
\address[A]{Department of Applied Mathematics\\
Hong Kong Polytechnic University\\
Kowloon \\
Hong Kong\\
\printead{e1}}

\address[B]{Mathematical Institute
and Oxford--Man\\
\quad Institute of Quantitative Finance\\
University of Oxford\\
24--29 St Giles\\
Oxford, OX1 3LB\\
United Kingdom\\
and\\
Department of Systems Engineering\\
\quad and Engineering Management\\
Chinese University of Hong Kong \\
Shatin \\
Hong Kong\\
\printead{e2}}

\end{aug}

\received{\smonth{3} \syear{2011}}
\revised{\smonth{12} \syear{2011}}

%
\begin{abstract}
We formulate an optimal stopping problem for a geometric Brownian
motion where the probability scale is distorted by a general nonlinear
function. The problem is inherently time inconsistent due to the
Choquet integration involved. We develop a new approach, based on a
reformulation of the problem where one optimally chooses the
probability distribution or quantile function of the stopped state. An
optimal stopping time can then be recovered from the obtained
distribution/quantile function, either in a straightforward way for
several important cases or in general via the Skorokhod embedding. This
approach enables us to solve the problem in a fairly general manner
with different shapes of the payoff and probability distortion functions.
We also discuss economical interpretations of the results. In
particular, we justify several liquidation strategies widely adopted in
stock trading,
including those of ``buy and hold,'' ``cut loss or take profit,'' ``cut
loss and let profit run'' and ``sell on a percentage of historical high.''
\end{abstract}




\end{frontmatter}

\section{Introduction}\label{sec1}

Many experimental evidences show that people tend to inflate,
intentionally or unintentionally, small probabilities.
Here we present two simplified examples.
We write a random variable (\textit{prospect})
$X=(x_i,p_i;\break i=1,2,\ldots,m)$ if $X=x_i$ with probability $p_i$, and
write $X\succ Y$ if prospect $X$ is preferred than prospect $Y$.
Then it is a general observation that $(\$5000,0.001; \$0,0.999)\succ(\$
5,1)$ although the two prospects have the same mean. One of the
explanations is that people usually exaggerate the small probability
associated\vadjust{\goodbreak} with a big payoff (so people buy lotteries). On the other
hand, it is common that $(-\$5,1)\succ(-\$5000,0.001;\$0,0.999)$,
indicating an inflation of the small probability in respect of a big
loss (so people buy insurances).

Probability distortion (or weighting) is one of the building blocks of
a number of modern behavioral economics theories including Kahneman and
Tversky's cumulative prospect theory (CPT) \cite
{Kahneman&Tversky,T&K1992} and Lopes' SP/A theory~\cite{Lopes}.
Yaari's dual theory of choice~\cite{Yaari} uses probability distortion
as a substitute for
expected utility in describing people's risk preferences. Probability
distortion has also been extensively investigated in the insurance
literature (see, e.g.,~\cite{Wang,WY,CMM}).

In this paper we introduce and study optimal stopping of a geometric
Brownian motion when probability scale is distorted. To our best
knowledge such a problem has not been formally formulated nor attacked
before. Due to the probability distortion, the payoff functional of the
stopping problem is evaluated via
the so-called Choquet integration, a type of nonlinear expectation. We
are interested in developing a general approach to solving the problem
and in
understanding whether and how the probability distortion changes
optimal stopping strategies.

There have been well-developed approaches in solving classical optimal
stopping without probability distortion, including those of probability
(martingale) and PDE (dynamic programming or variational inequality).
We refer to~\cite{Fri} and~\cite{Shr} for classical accounts of the theory.
These approaches are based crucially on the time-consistency of the
underlying problem. In~\cite{NO} and~\cite{Riedel}, the authors
study optimal stopping problems under Knightian uncertainty (or
ambiguity), involving essentially a different type of nonlinear
expectation in their payoff functionals.
However, both papers \textit{assume} upfront that time-consistency
(or, equivalently, the so-called rectangularity) is kept intact, which
enables the applicability of the classical approaches. Henderson \cite
{Vicky} investigates the
disposition effect in stock selling through an optimal stopping with
$S$-shaped payoff functions (motivated by Kahneman and Tversky's CPT);
however, since there is no probability distortion involved she is again
able to apply the martingale theory to solve the problem.\looseness=-1

In the presence of probability distortion, however,
the fundamental time-consistency structure is lost, to which the
traditional martingale or dynamic programming approaches fail to apply.
This is the major challenge arising from probability distortion in
optimal stopping.
Barberis~\cite{B} studies optimal exit strategies in casino gambling
with CPT preferences (including probability distortion) in a
discrete-time setting. He highlights the inherent time-inconsistency
issue of the problem and obtains only numerical solutions via
exhaustive enumeration.

In this paper we develop a new approach to overcome the difficulties
resulting from the (probability) distortion including the time
inconsistency. An important\vadjust{\goodbreak} technical ingredient in our approach is the
\textit{Skorokhod embedding}. Skorokhod~\cite{Skorokhod} introduced
and solved the following problem:
Given a standard Brownian motion~$B_t$ and a probability measure $m$
with 0 mean and
finite second moment, find an integrable stopping time $\tau$
such that the distribution of~$B_{\tau}$ is $m$. Since then there have been
great number of variants, generalizations and applications of the
Skorokhod embedding problem (see~\cite{Jan} for a recent survey).

Suppose the stochastic process to be stopped is $\{S_t,t\geqslant 0\}
$. The key idea in solving our ``distorted'' optimal stopping consists
of first determining the probability distribution of an optimally
stopped state,
$S_{\tau^*}$, and then recovering an optimal stopping $\tau^*$,
either in an obvious way for several important
cases or in general via the Skorokhod embedding.
The first part is inspired by the observation that
the payoff functional, even though evaluated under the distorted
probability, still depends only on the distribution function of the
stopped state $S_{\tau}$; so one can take the distribution function
(instead of the stopping time) as the decision variable in solving the
optimal stopping problem. The resulting problem is said to have a
\textit{distribution formulation}. In some cases
it is more convenient to consider the quantile function (the
left-continuous inverse of the distribution function) as the decision
variable, based on which we have the \textit{quantile
formulation}.\setcounter{footnote}{3}\footnote{Quantile formulation has been introduced and
developed in the context of financial portfolio
selection involving probability distortion (see~\cite{Schied1,DANA}
and~\cite{CD}
for earlier works). Jin and Zhou~\cite{Jin&Zhou} employ the
formulation to solve a continuous-time portfolio selection model with
the behavioral CPT preferences.
The quantile formulation has recently been further developed in \cite
{HZ0} into
a general paradigm of solving nonexpected utility maximization models.}
To summarize, our original problem can be generally solved by a
three-step procedure. The first step is to rewrite the problem in a
distribution or quantile formulation, the second one is to solve the
resulting distribution/quantile optimization problem and the last one
is to derive an optimal stopping from the
optimal distribution/quantile function.

The remainder of the paper is organized as follows. In Section~\ref{sec2}, we
formulate the optimal stopping problem under probability distortion,
and then transfer the problem into one where the underlying process is
a martingale.
In Section~\ref{sec3} we present the distribution and quantile formulations of
the original problem. In Sections~\ref{seconvex}--\ref{sec6} we solve the problem (resp., for
different shapes)\footnote{Throughout this paper the term ``shape''
mainly refers to the property of a function related to piecewise
convexity and concavity. A function is called
$S$-shaped (resp., reverse $S$-shaped) if it includes two pieces, with
the left piece being convex (resp., concave) and the right one concave
(resp., convex). These shapes all have economical interpretations
related to risk preferences.} of the probability distortion and the
payoff functions.
We also discuss financial/economical implications of the derived
results and compare our results with the case when there is no\vadjust{\goodbreak}
probability distortion.
In particular, we justify several liquidation strategies widely adopted
in stock trading.
We finally conclude this paper in Section~\ref{sec7}. Some technical proofs are
placed in Appendices~\ref{AA}--\ref{AE}.

\section{Optimal stopping formulation}\label{sec2}
\subsection{The problem}\label{sec2.1}
Consider a stochastic process, $\{P_t,t\geqslant 0\}$, that follows a~geometric Brownian motion (GBM)
%
\begin{equation}\label{gbm}
\dd P_{t}=\mu P_{t}\,\dd t+\sigma P_{t}\,\dd B_{t},\qquad P_0>0,
\end{equation}
where $\mu$ and $\sigma>0$ are real constants, and $\{B_t,t\geqslant
0\}$ is a standard one-dimensional Brownian motion in a filtered
probability space $(\Omega, \mathcal{F},\{\mathcal{F}_{t}\}
_{t\geqslant 0},\BP)$. In many discussions below, $\{P_t,t\geqslant
0\}$ will be interpreted as the price process of an asset.

Let $\mathcal{T}$ be the set of all $\{\mathcal{F}_{t}\}_{t\geqslant
0}$-stopping times $\tau$ with $\BP(\tau<+\infty)=1$.
A~decision-maker (agent) chooses $\tau\in\mathcal{T}$ to stop the
process and obtain a payoff $U(P_\tau)$, where $U(\cdot)\dvtx\R
^{+}\mapsto\R^{+}$ is a given nondecreasing, continuous function.
The agent distorts the probability scale with a distortion (weighting)
function $w(\cdot)\dvtx[0,1]\mapsto[0,1]$, which is a
strictly increasing, absolutely continuous function with $w(0)=0$ and $w(1)=1$.
The agent's target is to maximize her ``distorted'' mean payoff functional:
%
\begin{equation} \label{objective0}
\mbox{Maximize}\quad J(\tau):=\int_{0}^{\infty}w\bigl(\BP\bigl(U(P_{\tau})>
x\bigr)\bigr)\,\dd x
\end{equation}
over $\tau\in\mathcal{T}$. In probabilistic terms the above
criterion \eqref{objective0} is a nonlinear expectation, called the
\textit{Choquet expectation} or \textit{Choquet integral}, of the
random payoff $U(P_{\tau})$ under the \textit{capacity} $w(P(\cdot
))$. Note that when $U(P_{\tau})$ is a discrete random variable,
\eqref{objective0} agrees with that in the CPT~\cite{T&K1992} so our
criterion is a natural generalization of the CPT value function
covering both
continuous and discrete payoffs.

Another important point to note is that here the underlying process
$P_t$ is independent of the probability distortion. In the context of
stock trading,
this means that the agent is a ``small investor'' so her preference
only affects {her} own stopping strategies,
but not the asset dynamics. How probability distortions of market
participants might collectively affect the asset price
is a significant open problem and is certainly beyond the scope of this
paper.

If there is no probability distortion [i.e., $w(x)\equiv x$], then the
objective functional~\eqref{objective0} is nothing else than the
expected payoff appearing in a standard optimal stopping problem
\[
J(\tau)=\int_{0}^{\infty}w\bigl(\BP\bigl(U(P_{\tau})> x\bigr)\bigr)\,\dd x=\int
_{0}^{\infty}\BP\bigl(U(P_{\tau})> x\bigr)\,\dd x=\BE[U(P_{\tau})].
\]
Hence, the problem considered in this paper is that of a ``distorted''
optimal stopping in the sense that\vadjust{\goodbreak} the probability scale is distorted.
As with the classical optimal stopping there can be many applications
of our formulation. For instance, the following problem falls into our
formulation:
An agent with CPT preferences needs to determine the time of exercising
a perpetual American option written on an
asset whose discounted price process follows \eqref{gbm}, whereas the
option pays $U(P)$ at the exercise price $P$.\footnote{We assume in
this paper that $U(\cdot)$ is nondecreasing. Although some payoff
functions may be nonincreasing, such as that of a put option,
the case of a nonincreasing $U(\cdot)$ can be dealt with in exactly
the same way as with the nondecreasing counterpart to be presented in
this paper. On the other hand, we do not assume
$U(\cdot)$ to be smooth or strictly increasing so as to accommodate
call option type of payoffs.}
Our problem can also be interpreted simply as an investor hoping to
determine the best selling time of a stock that she is holding, and
$U(\cdot)$ in this case is a
utility function of the proceeds of the liquidation. Yet another
example of our formulation is the so-called irreversible investment
where the objective is to determine the best time to carry out an
investment project (see, e.g.,~\cite{DP,NO}).

\subsection{Transformation}\label{sec2.2}
For subsequent analysis we need to transform problem~\eqref
{objective0} into one where the underlying\vspace*{1pt} process is a martingale.
To proceed let us first study the simple case when $\mu=\frac
{1}{2}\sigma^{2}$. Indeed, in this case $P_{t}=P_{0}e^{\sigma B_{t}}$.
Let
\[
\tau_{x}=\inf\{t\geqslant 0\dvtx B_{t}= \sigma^{-1}\ln(x/P_{0})\}
\qquad\forall x\in(0,+\infty).
\]
Then $\tau_{x}\in\mathcal{T}$, $P_{\tau_{x}}=x$ almost surely, and
$J(\tau_{x})=U(x)$, $\forall x\in(0,+\infty)$. However, for any
$\tau\in\mathcal{T}$,
\begin{eqnarray*}
J(\tau)&=&\int_{0}^{\infty}w\bigl(\BP\bigl(U(P_{\tau})> x\bigr)\bigr)\,\dd x=\int
_{0}^{\bar U}w\bigl(\BP\bigl(U(P_{\tau})> x\bigr)\bigr)\,\dd x\\
&\leqslant&\int_{0}^{\bar U}w(1)\,\dd x= \bar U =\sup_{x>0} J(\tau_{x}),
\end{eqnarray*}
where $\bar U:=\sup_{x>0} U(x)$.
This shows that the optimal value of problem \eqref{objective0} is~$\bar U$,
and that an optimal stopping time, if it ever exists, is of the form
$\tau_x$.
Moreover, if there exists at least one $x^{*}>0$ such that $U(x^{*})=
\bar U$, then $\tau_{x^{*}}$ is an optimal stopping time. If, on the
other hand, $U(y)<\bar U$ for every $y>0$,\vspace*{1pt}
then for any stopping time $\tau\in\mathcal{T}$, we have $U(P_{\tau
})<\bar U$. Therefore, noting that $w(\cdot)$ is strictly increasing,
\begin{eqnarray*}
J(\tau)=\int_{0}^{\infty}w\bigl(\BP\bigl(U(P_{\tau})> x\bigr)\bigr)\,\dd x
<\int_{0}^{\infty}w\bigl(\BP(\bar U> x)\bigr)\,\dd x
= \bar U,
\end{eqnarray*}
which means that the optimal value is not achievable by any stopping
time. However, $\lim_{n\to\infty}J(\tau_{x_n})=\sup_{\tau\in
\mathcal{T}}J(\tau)$ for any
sequence $\{x_n>0\dvtx n=1,2,\ldots\}$\vspace*{1pt} satisfying $\lim_{n\to\infty
}U(x_n)=\bar U$.

Given that the case when $\mu=\frac{1}{2}\sigma^{2}$ has been
completely solved, henceforth we only consider the case when $\mu\neq
\frac{1}{2}\sigma^{2}$. We now\vadjust{\goodbreak} convert problem \eqref{objective0}
into an equivalent one.
Let
%
\begin{eqnarray}\label{defu}
\beta&:=&\frac{-2\mu+\sigma^{2}}{\sigma^{2}}\neq0,\qquad
S_{t}:=P_{t}^{\beta},
\nonumber
\\[-9pt]
\\[-9pt]
\nonumber
 u(x)&:=&U(x^{1/\beta})\qquad \forall x\in(0,+\infty).
\end{eqnarray}
Then It\^o's rule gives
%
\begin{equation}\label{eqs}
\dd S_{t}= \beta\sigma S_{t}\,\dd B_{t},\qquad S_{0}=P_{0}^{\beta}:=s>0.
\end{equation}
Now we can rewrite problem \eqref{objective0} as
%
\begin{eqnarray}\label{objective}
\qquad\mbox{Maximize}\quad J(\tau)&=&\int_{0}^{\infty}w\bigl(\BP\bigl(U(P_{\tau})>
x\bigr)\bigr)\,\dd x
\nonumber
\\[-9pt]
\\[-9pt]
\nonumber
&=&\int_{0}^{\infty}w\bigl(\BP(u(S_{\tau})> x)\bigr)\,\dd x
\end{eqnarray}
over $\tau\in\mathcal{T}$,
where the new, auxiliary process $S_t$ follows \eqref{eqs}, and the
new payoff function $u(\cdot)$ is defined in \eqref{defu}. In the
remainder of this paper we will mainly consider the objective
functional in \eqref{objective} instead of that in \eqref{objective0}.

The advantage of this transformation is that $S_t$ is now a martingale,
which enables us to apply the Skorokhod theorem later on.
Interestingly, $u(\cdot)$ may~now have a completely different shape
than $U(\cdot)$, depending on the value of~$\beta$.

\subsection{Examples}\label{sec2.3}
We now discuss several popular payoff functions $U(\cdot)$ as
examples; these examples will also serve as benchmarks for illustrating
the main results of this paper.

Let us start with the example of a call option written on an underlying
asset whose discounted price process follows (\ref{gbm}). The payoff function
is $U(x)=(x-K)^{+}$ for some $K>0$. Then $u(x)=(x^{1/\beta}-K)^{+}$.
If $ \beta<0 $ or equivalently $\frac{\mu}{\sigma^2}>0.5$, the
underlying asset is ``good.''\footnote{In~\cite{SXZ}, $\frac{\mu
}{\sigma^2}$ is termed the ``goodness index'' of an asset.} In this
case $u(\cdot)$ is nonincreasing and convex. If
$0<\beta\leqslant1$ or $0\leqslant\frac{\mu}{\sigma^2}<0.5$, the
asset is between good and bad and $u(\cdot)$ is nondecreasing and convex.
If $\beta>1$ or $\mu<0$, the asset is ``bad,'' and $u(\cdot)$ is
nondecreasing and $S$-shaped.

Now take a power function $U(x)=\frac{1}{\gamma}x^{\gamma}$, $\gamma
\in(0,1)$. Then $u(x)=\frac{1}{\gamma}x^{\gamma/\beta}$ is
strictly decreasing and convex if $\beta< 0$, strictly increasing and
convex if $0<\beta\leqslant\gamma$ (i.e., the asset is ``not so
bad'' in respect of the original payoff/utility function) and strictly
increasing and concave
if $\beta>\gamma$ (the asset is sufficiently bad).

For a log utility function $U(x)=\ln(x+1)$, $u(x)=\ln(x^{1/\beta
}+1)$ is strictly decreasing if $\beta < 0$, strictly increasing and
$S$-shaped if $0<\beta<1$ and strictly increasing and concave if
$\beta\geqslant 1$.\vadjust{\goodbreak}

For an exponential utility function $U(x)=1-e^{-\alpha x}$, $\alpha
>0$, $u(x)=1-e^{-\alpha x^{1/\beta}}$ is strictly decreasing if
$\beta < 0$, strictly increasing and $S$-shaped if $0<\beta< 1$ and
strictly increasing and concave if $\beta\geqslant 1$.

Next, let us take an $S$-shaped piecewise power function
$U(x)=(x/k)^{\alpha_{1}}\times \mathbf{1}_{(0,k]}(x)+(x/k)^{\alpha
_{2}}\mathbf{1}_{(k,\infty)}(x)$, where\vspace*{1pt} $ \alpha_{1}\geqslant
1\geqslant \alpha_{2}>0$, $k>0$.
If $\beta < 0$, then $u(x)=x^{\alpha_{2}/
\beta}k^{-\alpha_{2}}\mathbf{1}_{(0,k^{\beta}]}(x)+x^{\alpha_{1}/\beta}k^{-\alpha_{1}}\mathbf{1}_{(k^{\beta},\infty)}(x)$
is strictly decreasing.
Otherwise,  when $\beta\geq0$, $u(x)=x^{\alpha_{1}/\beta}k^{-\alpha_{1}}\mathbf{1}_{(0,k^{\beta}]}(x)
+x^{\alpha_{2}/\beta}k^{-\alpha_{2}}\mathbf{1}_{(k^{\beta},\infty)}(x)$ is strictly increasing and piecewise convex if $0<\beta <\alpha_{2}$,
strictly increasing and  $S$-shaped if $\alpha_{2}\leq \beta \leq \alpha_{1}$, and strictly increasing and  concave
if $\beta > \alpha_{1}$.

Finally, for a
general nondecreasing function $U(\cdot)$, $u(x)=U(x^{1/\beta})$ is
nonincreasing if $\beta< 0$ and nondecreasing if $\beta> 0$.

\subsection{Solution to a trivial case}\label{sec2.4}
While solving \eqref{objective} in general requires a new approach,
which will be developed in
the subsequent sections, in this subsection we present the solution to
a mathematically (almost) trivial yet economically significant case.\vspace*{-1pt}

\begin{thmm}\label{bah}
If $u(\cdot)$ is nonincreasing, then problem \eqref{objective} has
the optimal value $u(0+)$ and
%
\begin{equation}\label{infinite}
\lim_{T\to+\infty}J(T)=\sup_{\tau\in\mathcal{T}}J(\tau).
\end{equation}
Moreover, if $u(\ell)=u(0+)$ for some $\ell>0$, then
$\tau_{\ell}:=\inf\{t\geqslant 0\dvtx  S_{t}\leqslant\ell\}$
is an optimal stopping time for problem \eqref{objective}. If $u(\ell
)<u(0+)$ for every $\ell>0$, then~\eqref{objective} has no optimal solution.
\end{thmm}

A proof can be found in Appendix~\ref{AA}. We remark that $u(\cdot)$
is not required to be even continuous in the proof.

Identity \eqref{infinite} suggests that the supremum of the payoff
functional can be achieved by \textit{not} stopping at all, if
$u(\cdot)$ is nonincreasing.
There is an interesting economical interpretation of the above result
in the context of asset selling. In all the examples presented in
Section~\ref{sec2.3}, the case of $u(\cdot)$ being nonincreasing corresponds to
$\beta<0$, namely, the underlying asset being good. Moreover, in all
but the last general example, it holds that
$u(0+)>u(\ell)$ for all $\ell>0$. Theorem~\ref{bah} then indicates
that one should not sell at any price level, or one should hold the
asset perpetually.
This is indeed consistent with the traditional investment wisdom that
one should ``buy and hold a good asset.''\footnote{In~\cite{SXZ},
a similar result is derived, albeit for a different asset selling model
where the time horizon is finite, probability distortion is absent and
the objective is to minimize the relative error between the selling
price and
the all-time-high price.}

\section{Distribution/quantile formulation}\label{sec3}
In view of Theorem~\ref{bah}, henceforth we consider only the case
when $u(\cdot)$ is nondecreasing. Let us specify the standing
assumption we impose from this point on.

\begin{assmp}
$u(\cdot)\dvtx  \R^{+}\mapsto\R^{+}$ is nondecreasing, absolutely
continuous with $u(0)=0$; $w(\cdot)\dvtx [0,1]\mapsto[0,1]$ is
strictly increasing, absolutely continuous with $w(0)=0$ and $w(1)=1$.
\end{assmp}

Note that $u(0)=0$ is just for simplicity, as one may consider $\bar
{u}(\cdot)=u(\cdot)-u(0)$ if $u(0)\neq0$.

Throughout this paper, for any nondecreasing function $f\dvtx \R^{+}\mapsto
[0,1]$, we denote by $f^{-1}\dvtx [0,1)\mapsto\R^{+}$ the left-continuous
inverse function of $f$ which is defined by
\[
f^{-1}(x):=\inf\{y\in\R^{+}\dvtx  f(y)\geqslant x\},\qquad x\in[0,1).
\]
Clearly $f^{-1}$ is nondecreasing and left-continuous.
We say $F\dvtx \R^{+}\mapsto[0,1]$ is a cumulative distribution function
(CDF) if $F(0)=0, F(+\infty)\equiv\lim_{x\rightarrow+\infty
}F(x)=1$ and $F$ is nondecreasing and \textit{c\'adl\'ag}.
We call $G\dvtx  [0,1)\mapsto\R^{+}$ a quantile function if $G(0)=0$,
$G(x)>0\ \forall x\in(0,1)$, $G$ is nondecreasing and
left-continuous.\footnote{Note that in this paper the underlying
process $S_t$ is strictly positive at any time; hence,
we need to consider only the CDF and quantile function of strictly
positive random variables.}

Now define the following distribution set $\mathcal{D}$ and quantile
set $\mathcal{Q}$:
\begin{eqnarray*}
\mathcal{D}&:=&\{F\dvtx \R^{+}\mapsto[0,1] | F\ \mathrm{is\ the\ CDF\
of\ }S_{\tau},\ \mathrm{for\ some\ }\tau\in\mathcal{T}\},\\
\mathcal{Q}&:=&\{G\dvtx [0,1)\mapsto\R^{+} | G=F^{-1}\ \mathrm{for\
some\ } F\in\mathcal{D}\}.
\end{eqnarray*}

%
\begin{lemma} \label{2forms}
For any $\tau\in\mathcal{T}$, we have
%
\begin{eqnarray}
\label{dis-form} J(\tau)&=&J_{D}(F):=\int_{0}^{\infty
}w\bigl(1-F(x)\bigr)u'(x)\,\dd x,\\
\label{qua-form} J(\tau)&=&J_{Q}(G):=\int_{0}^{1}u(G(x))w'(1-x)\,\dd x,
\end{eqnarray}
where $F$ and $G$ are the CDF and the quantile function of $S_{\tau}$,
respectively. Moreover,
%
\begin{equation}\label{equiv}
\sup_{\tau\in\mathcal{T}}J(\tau)=\sup_{F\in\mathcal
{D}}J_{D}(F)=\sup_{G\in\mathcal{Q}}J_{Q}(G).
\end{equation}
\end{lemma}


A proof is relegated to Appendix~\ref{AB}.

We name \eqref{dis-form} and \eqref{qua-form} as the \textit
{distribution formulation} and the \textit{quantile formulation} of
problem \eqref{objective}, respectively.\vadjust{\goodbreak}

We observe certain symmetry (or rather duality) between the
distribution formulation \eqref{dis-form} and the quantile formulation
\eqref{qua-form}. In particular, $w(\cdot)$ and $u(\cdot)$ play
symmetric roles in the two formulations.\footnote{Indeed, in the
context of utility theory both a probability
distortion function and a utility function describe an investor's
preference toward risk; they do play some dual roles (see~\cite{Yaari}).}
The availability of two formulations enables us to choose a convenient
one in solving the original stopping problem~\eqref{objective},
depending on the
shape of $w(\cdot)$ and $u(\cdot)$. For instance, if $u(\cdot)$ is
known to be concave or convex [while $w(\cdot)$ is arbitrary], then
it might be advantageous
to choose the quantile formulation \eqref{qua-form}.

Next, we need to characterize the sets $\mathcal{D}$ and $\mathcal
{Q}$ more explicitly for the second step (solving the
distribution/quantile optimization problem).

Let $F\in\mathcal{D}$, namely, $F$ is the CDF of $S_\tau$ for some
$\tau\in\mathcal{T}$.
Since $S_t$ is a nonnegative martingale, optional sampling theorem and
Fatou's lemma yield, necessarily, $\int_{0}^{\infty}(1-F(x))\,\dd
x\equiv\BE[S_{\tau}]\leqslant s$.
It turns out that this inequality, $\int_{0}^{\infty}(1-F(x))\,\dd
x\leqslant s$,
is not only necessary but also sufficient for $F$ to belong to
$\mathcal{D}$.
%
\begin{lemma}
We have the following expressions of the distribution set $\mathcal
{D}$ and quantile set $\mathcal{Q}$:
\begin{eqnarray*}
\mathcal{D}&=&\biggl\{F\dvtx \R^{+}\mapsto[0,1] \Big|\mbox{$F$ is a CDF and } \int
_{0}^{\infty}\bigl(1-F(x)\bigr)\,\dd x\leqslant s\biggr\},\\[-2pt]
\mathcal{Q}&=&\biggl\{G\dvtx [0,1)\mapsto\R^{+} \Big|\mbox{$G$ is a quantile
function and } \int_{0}^{1}G(x)\,\dd x\leqslant s\biggr\}.
\end{eqnarray*}
\end{lemma}

\begin{pf}
First assume $\beta>0$. We write $S_{t}=s\exp(-\frac{1}{2}\beta
^2\sigma^{2}t+ \beta\sigma B_{t})\equiv s\exp(\beta\sigma\tilde
B_{t})$, where $\tilde B_{t}:=B_{t}-\frac{1}{2}\beta\sigma t$ is a
drifted Brownian motion with a \textit{negative} drift. Denote by
$F_X$ the CDF of a random variable $X$. For any $\tau\in\mathcal
{T}$, we have
\[
F_{\tilde B_{\tau}}(x)=\BP(\tilde B_{\tau}\leqslant x)=\BP(S_{\tau
}\leqslant se^{ \beta\sigma x})
=F_{S_{\tau}}(se^{\beta\sigma x}).
\]
On the other hand, according to Theorem 2.1 in~\cite{Hall}, a CDF $F$
is the CDF of $\tilde B_{\tau}$ for some $\tau\in\mathcal{T}$ if
and only if
$\int_{-\infty}^{\infty}e^{\beta\sigma x}\,\dd F(x)\leqslant1.$
So $F$ is the CDF of $S_{\tau}$ for some $\tau\in\mathcal{T}$ if
and only if it is a CDF and
$
\int_{-\infty}^{\infty}e^{\beta\sigma x}\,\dd F(se^{\beta\sigma
x})\leqslant1,
$
or
$
\int_{0}^{\infty}x\,\dd F(x)\leqslant s.
$
The above is equivalent to
$\int_{0}^{\infty}(1-F(x))\,\dd x\leqslant s $, or $\int
_{0}^{1}G(x)\,\dd x\leqslant s.$

Now if $\beta<0$, then write $S_{t}=s\exp(-\beta\sigma\tilde
B_{t})$, where $\tilde B_{t}:=-B_{t}+\frac{1}{2}\beta\sigma t$ is
still a drifted Brownian motion with a negative drift. The rest of the
proof is exactly the same as above.
This completes the proof.
\end{pf}

An important by-product of this lemma is that both the sets $\mathcal
{D}$ and $\mathcal{Q}$ are convex.\vadjust{\goodbreak}

Now we have reformulated our problem \eqref{objective} into
optimization problems maximizing \eqref{dis-form} or \eqref{qua-form}
over a convex set $\mathcal{D}$ or $\mathcal{Q}$, respectively. In
the following several sections, we will solve these problems with
different shapes of the functions $u(\cdot)$ and
$w(\cdot)$.

\section{\texorpdfstring{Convex $w(\cdot)$ or $u(\cdot)$}{Convex $w(.)$ or $u(.)$}}\label{seconvex}
In this section, we will solve problem \eqref{objective} assuming that
\textit{either} $w(\cdot)$ \textit{or} $u(\cdot)$ is convex.

\subsection{\texorpdfstring{Convex $w(\cdot)$}{Convex $w(.)$}}\label{sec4.1}
Assume for now that $w(\cdot)$ is convex [whereas $u(\cdot)$, being
nondecreasing, is allowed to have any shape].
In this case the distribution formulation \eqref{dis-form} is easier
to study than its quantile counterpart \eqref{qua-form},
since the shape of $u(\cdot)$ is unknown in the latter. The
distribution formulation in this case is to maximize a convex
functional over a convex set. Intuitively speaking, a
maximum of~\eqref{dis-form} should be at ``corners'' of the constraint
set, $\mathcal{D}$. We are going to establish that
these corners must be step functions having at most two jumps.

For $n=2,3,\ldots,$ define
\[
\mathcal{S}_{n}:=\Biggl\{F \dvtx F=\sum_{i=1}^{n-1} c_{i}\mathbf
{1}_{[a_{i},a_{i+1})}+\mathbf{1}_{[a_{n},\infty)}, 0< c_{i}\leqslant
c_{i+1} \leqslant1, 0 < a_{i}\leqslant a_{i+1} \Biggr\}
\]
and $\mathcal{D}_{n}:=\mathcal{S}_{n}\cap\mathcal{D}$. Clearly
$\mathcal{S}_{n}\subseteq\mathcal{S}_{n+1}$ and $\mathcal
{D}_{n}\subseteq\mathcal{D}_{n+1}\subseteq\mathcal{D} $,
$n=2,3,\ldots.$

\begin{lemma} \label{f=f2}
If $w(\cdot)$ is convex, then
\[
\sup_{F\in\mathcal{D}}J_{D}(F)=\sup_{F\in\mathcal{D}_{2}}J_{D}(F).
\]
\end{lemma}

A proof of this lemma is provided in Appendix~\ref{AC}.

By virtue of Lemma~\ref{f=f2}, in maximizing \eqref{dis-form} we need
only to search over the set $\mathcal{D}_{2}$ or, equivalently, to
find the best parameters $a$, $b$ and $c$ in defining an element
$F(x)=c\mathbf{1}_{[a,b)}(x)+\mathbf{1}_{[b,+\infty)}(x)$ in
$\mathcal{D}_{2}$. This becomes a three-dimensional constrained
optimization problem which is dramatically easier to solve than the
original stopping problem.
We present the results in the following theorem.

\begin{thmm}\label{convexw}
If $w(\cdot)$ is convex, then
\[
\sup_{\tau\in\mathcal{T}}J(\tau)=\sup_{0<a \leqslant s \leqslant
b}\biggl[\biggl(1-w\biggl(\frac{s-a}{b-a}\biggr)\biggr)u(a)+w\biggl(\frac{s-a}{b-a}\biggr)u(b)\biggr].
\]
Moreover, if
%
\begin{equation}\label{abstar}
(a^{*}, b^{*})=\mathop{\argmax}_{0< a \leqslant s \leqslant b}\biggl[\biggl(1-w\biggl(\frac
{s-a}{b-a}\biggr)\biggr)u(a)+w\biggl(\frac{s-a}{b-a}\biggr)u(b)\biggr],
\end{equation}
then
%
\begin{equation}\label{taustar}
\tau_{(a^{*}, b^{*})}:=\cases{
\inf\{t\geqslant 0\dvtx S_{t}\notin(a^{*}, b^{*})\},&\quad $\mbox{if }
a^{*}<b^{*}$,\vspace*{2pt}\cr
0,&\quad$\mbox{if } a^{*}=b^{*}$}
\end{equation}
is an optimal stopping to problem \eqref{objective}.
\end{thmm}
\begin{pf}
Due to Lemma~\ref{f=f2}, we need only to find the optimal distribution
function in $\mathcal{D}_{2}$ to maximize \eqref{dis-form}. For any
$F\in\mathcal{D}_{2}$ with
$F(x)=c\mathbf{1}_{[a,b)}(x)+\mathbf{1}_{[b,+\infty)}(x), x\in
[0,+\infty),$
we have
\[
J_{D}(F)=\int_{0}^{\infty}w\bigl(1-F(x)\bigr)u'(x)\,\dd x=\bigl(1-w(1-c)\bigr)u(a)+w(1-c)u(b)
\]
and
\[
\int_{0}^{\infty}\bigl(1-F(x)\bigr)\,\dd x=ac+b(1-c).
\]
Thus our problem boils down to
%
\begin{eqnarray}\label{lp}
&&\mbox{Maximize}\quad J(a,b,c):=\bigl(1-w(1-c)\bigr)u(a)+w(1-c)u(b)
\nonumber
\\[-8pt]
\\[-8pt]
\nonumber
&&\quad\mbox{subject to:}\quad ac+b(1-c)\leqslant s,\ 0< a \leqslant b,\
0\leqslant c\leqslant1.
\end{eqnarray}
Clearly $a\leqslant s$, otherwise the first constraint of \eqref{lp}
would be violated. On the other hand, in maximizing $J(a,b,c)$ one
should choose $b$ as large as possible when $a$ and $c$ are fixed. So
we need only to consider the range $0<a \leqslant s \leqslant b$ when solving
\eqref{lp}. Moreover, $J(a,b,c)$ is nonincreasing in $c$ when $a$ and
$b$ are fixed; hence, $c=\frac{b-s}{b-a}$ when $a<b$, while
$c\in[0,1]$ can be arbitrarily chosen
when $a=b$.
Therefore,
\[
\sup_{\tau\in\mathcal{T}}J(\tau)=\sup_{F\in\mathcal{D}_{2}}J_{D}(F)
=\sup_{0< a \leqslant s \leqslant b}\biggl[\biggl(1-w\biggl(\frac
{s-a}{b-a}\biggr)\biggr)u(a)+w\biggl(\frac{s-a}{b-a}\biggr)u(b)\biggr].
\]

Now if $(a^{*}, b^{*})$ with $0<a^{*}<b^*$ is determined by \eqref
{abstar}, then clearly $S_{\tau_{(a^{*}, b^{*})}}$, where
$\tau_{(a^{*}, b^{*})}$ is defined by \eqref{taustar}, has a
two-point distribution\break $P(S_{\tau_{(a^{*}, b^{*})}}=a^*)=c^*$ and
$P(S_{\tau_{(a^{*}, b^{*})}}=b^*)=1-c^*$. Moreover, $c^*=\frac
{b^*-s}{b^*-a^*}$ by virtue of the optional sampling theorem. The CDF
of $S_{\tau_{(a^{*}, b^{*})}}$
is
$F^*(x)=c^*\mathbf{1}_{[a^*,b^*)}(x)+\mathbf{1}_{[b^*,+\infty)}(x)$.
Hence, $\tau_{(a^{*}, b^{*})}$ is an optimal solution. If
$a^{*}=b^{*}$, then it must hold that \mbox{$a^{*}=b^{*}=s$}. Hence, we have
$\sup_{\tau\in\mathcal{T}}J(\tau)=J(a^{*},b^{*},c)=u(s)=J(\tau
_{(a^{*}, b^{*})})$. So $\tau_{(a^{*}, b^{*})}$ defined by \eqref
{taustar} is an optimal solution to problem \eqref{objective}.
\end{pf}

According to~\cite{Yaari}, a convex probability distortion overweighs
``bad'' outcomes and underweighs ``good'' ones in maximizing the underlying
criterion; hence, it captures the risk-aversion of an investor.
The preceding theorem suggests that a risk averse agent's optimal
strategy is to stop at one of the two thresholds, $a^*$ and $b^*$. In
the context of stock liquidation, this
corresponds to the widely adopted\vadjust{\goodbreak} ``take-profit-or-cut-loss'' strategy,
namely, one should sell a stock either when it has reached a
pre-determined target $b^*$
or sunk to a prescribed loss level $a^*$ (note that the initial price
$s$ is in between $a^*$ and $b^*$), for a stock that is not worth ``buy
and hold perpetually.''

In particular, when there is no probability distortion [i.e.,
$w(x)\equiv x$] which is trivially convex, Theorem~\ref{convexw}
recovers the results of~\cite{Vicky} where an optimal stopping problem
is studied with a specific $S$-shaped utility function $u(\cdot)$
without probability distortion. Indeed,
Theorem~\ref{convexw} leads to a very general result
in the absence of
probability distortion: the optimality of the take-profit-or-cut-loss
strategy is inherent regardless of the shape of $u(\cdot)$ (be it
concave, convex or $S$-shaped) so long as it is nondecreasing.

It should be noted that in the current case no Skorokhod embedding
technique is needed to recover
the optimal stopping time $\tau^*$ from $S_{\tau^*}$. This is because
the explicit form of CDF of $S_{\tau^*}$ obtained reveals
that $S_{\tau^*}$ is a two-point distribution; hence, $\tau^*$ must
be the exit time of an interval.\vspace*{-2pt}

\begin{coro}\label{u-concave,w-convex}
If $u(\cdot)$ is concave and $w(\cdot)$ is convex, then
$
\sup_{\tau\in\mathcal{T}}J(\tau)=u(s).
$
Moreover, $\tau\equiv0$ is an optimal stopping.\vspace*{-2pt}
\end{coro}

\begin{pf}
The convexity of $w(\cdot)$ along with $w(0)=0$ and $w(1)=1$ implies
that $w(x)\leqslant x$, for all $x\in[0,1]$;
so we have
\begin{eqnarray*}
\sup_{\tau\in\mathcal{T}}J(\tau)&=&\sup_{0<a \leqslant s
\leqslant b}\biggl[\biggl(1-w\biggl(\frac{s-a}{b-a}\biggr)\biggr)u(a)+w\biggl(\frac{s-a}{b-a}\biggr)u(b)\biggr]\\
&=&\sup_{0<a \leqslant s \leqslant b}\biggl[u(a)+w\biggl(\frac
{s-a}{b-a}\biggr)\bigl(u(b)-u(a)\bigr)\biggr]\\
&\leqslant&\sup_{0<a \leqslant s \leqslant b}\biggl[u(a)+\frac
{s-a}{b-a}\bigl(u(b)-u(a)\bigr)\biggr]\\
&=&\sup_{0<a \leqslant s \leqslant b}\biggl[\frac{b-s}{b-a}u(a)+\frac
{s-a}{b-a}u(b)\biggr]
\leqslant u(s)=J(0),
\end{eqnarray*}
where we used the concavity property of $u(\cdot)$ to obtain the last
inequality.\vspace*{-2pt}~%
\end{pf}

This result stipulates that when $w(\cdot)$ is convex and $u(\cdot)$
is concave, the two thresholds degenerate into one which is the initial
state $s$. From some of the examples
in Section~\ref{sec2} (e.g., when the original payoff function is power or
logarithmic), $u(\cdot)$ being concave corresponds to a ``bad'' asset.
So if the agent is
risk averse [reflected by the convexity of $w(\cdot)$] or risk neutral
(no distortion) and the asset is unfavorable, then the optimal stopping
strategy is to stop immediately.

\subsection{\texorpdfstring{Convex $u(\cdot)$}{Convex $u(.)$}}\label{sec4.2}
Next we consider the case when $u(\cdot)$ is convex while $w(\cdot)$
has an arbitrary shape.
In this case the quantile formulation \eqref{qua-form} is more
convenient to deal with. The following result is an analog of Lemma~\ref{f=f2},
whose proof is, however, much simpler.\vadjust{\goodbreak}

\begin{lemma} \label{g=g2}
If $u(\cdot)$ is convex, then
\[
\sup_{G\in\mathcal{Q}}J_{Q}(G)=\sup_{G\in\mathcal{Q}_{2}}J_{Q}(G),
\]
where $\mathcal{Q}_{2}$ is defined as
\[
\mathcal{Q}_{2}:=\bigl\{G\in\mathcal{Q}\dvtx G=a\mathbf{1}_{(0,c]}+b\mathbf
{1}_{(c,1)}, 0<a\leqslant b, 0< c\leqslant1\bigr\}.
\]
\end{lemma}

See Appendix~\ref{AD} for a proof of the above lemma.

\begin{thmm}\label{uconvex}
If $u(\cdot)$ is convex, then
%
\begin{eqnarray}\label{g2}
\sup_{\tau\in\mathcal{T}}J(\tau)&=&\sup_{0< a \leqslant s \leqslant
b}\biggl[\biggl(1-w\biggl(\frac{s-a}{b-a}\biggr)\biggr)u(a)+w\biggl(\frac{s-a}{b-a}\biggr)u(b)\biggr]
\nonumber
\\[-8pt]
\\[-8pt]
\nonumber
&=&\sup_{x\in(0,1]}\biggl[ w(x)u\biggl(\frac{s}{x}\biggr)\biggr].
\end{eqnarray}
Moreover, if $
(a^{*}, b^{*})$ is determined by \eqref{abstar} then
$
\tau_{(a^{*}, b^{*})}$ defined by \eqref{taustar}
is an optimal stopping to problem \eqref{objective}.
\end{thmm}
\begin{pf}
Due to Lemma~\ref{g=g2}, we need only to find the optimal quantile
function in $\mathcal{Q}_{2}$ to maximize \eqref{qua-form}. For any
$G\in\mathcal{Q}_{2}$ with
$G(x)=a\mathbf{1}_{(0,c]}(x)+b\mathbf{1}_{(c,1)}(x),$ $x\in[0,1),$
we have
\begin{eqnarray*}
J_{Q}(G)=\int_{0}^{1}u(G(x))w'(1-x)\,\dd x=\bigl(1-w(1-c)\bigr)u(a)+w(1-c)u(b)
\end{eqnarray*}
and
\[
\int_{0}^{1}G(x)\,\dd x=ac+b(1-c).
\]
This leads to exactly the same optimization problem \eqref{lp}, and
one follows exactly the same lines of proof of Theorem~\ref{convexw}
to conclude that the first equality of \eqref{g2} is valid and
$\tau_{(a^{*}, b^{*})}$ defined by \eqref{taustar} is an optimal solution.

It remains to prove the second equality of \eqref{g2}.
Since both $u(\cdot)$ and $w(\cdot)$ are continuous, we have
\begin{eqnarray*}
&&\sup_{0< a \leqslant s \leqslant b}\biggl[\biggl(1-w\biggl(\frac
{s-a}{b-a}\biggr)\biggr)u(a)+w\biggl(\frac{s-a}{b-a}\biggr)u(b)\biggr]\\
&&\qquad\geqslant \sup_{a =0, s \leqslant b}\biggl[\biggl(1-w\biggl(\frac
{s-a}{b-a}\biggr)\biggr)u(a)+w\biggl(\frac{s-a}{b-a}\biggr)u(b)\biggr]
\\
&&\qquad=\sup_{x\in(0,1]}\biggl[ w(x)u\biggl(\frac{s}{x}\biggr)\biggr].
\end{eqnarray*}
Fix $0< a\leqslant s\leqslant b$ with $a<b$, and let
$G(x)=a\mathbf{1}_{(0,c]}(x)+b\mathbf{1}_{(c,1)}(x),$ $x\in[0,1),$
where $c=\frac{b-s}{b-a}$. Rewriting
\[
G(x)=\frac{a}{s}s+\frac{(b-a)(1-c)}{s}\frac{s}{1-c}\mathbf
{1}_{(c,1)}(x),\qquad x\in[0,1),
\]
we deduce by the convexity of $u(\cdot)$ that
\begin{eqnarray}
\nonumber u(G(x))&\leqslant&\frac{a}{s}u(s)+\frac
{(b-a)(1-c)}{s}u\biggl(\frac{s}{1-c}\mathbf{1}_{(c,1)}(x)\biggr),\qquad x\in[0,1);
\end{eqnarray}
hence,
%
\begin{eqnarray}
\label{g2<g1}
\qquad\int_{0}^{1}u(G(x))w'(1-x)\,\dd x&\leqslant&\frac
{a}{s}u(s)+\frac{(b-a)(1-c)}{s}u\biggl(\frac{s}{1-c}\biggr)w(1-c)
\nonumber
\\[-8pt]
\\[-8pt]
\nonumber
&\leqslant&\sup_{x\in(0,1]}\biggl[w(x)u\biggl(\frac{s}{x}\biggr)\biggr].
\end{eqnarray}
In other words,
\begin{eqnarray*}
\biggl(1-w\biggl(\frac{s-a}{b-a}\biggr)\biggr)u(a)+w\biggl(\frac{s-a}{b-a}\biggr)u(b)&=&\int
_{0}^{1}u(G(x))w'(1-x)\,\dd x\\
&\leqslant&\sup_{x\in(0,1]}\biggl[ w(x)u\biggl(\frac{s}{x}\biggr)\biggr].
\end{eqnarray*}
The proof is complete.
\end{pf}
%
\begin{coro}\label{care}
If $u(\cdot)$ is convex, then $\tau^{*}\equiv0$ is an optimal
solution to problem \eqref{objective} if and only if $u(s)=\sup_{x\in
(0,1]}[ w(x)u(\frac{s}{x})]$. Moreover, if $u(s)<\sup_{x\in(0,1]}[
w(x)u(\frac{s}{x})]$, then the maximum in \eqref{abstar} is not achievable.
\end{coro}

\begin{pf}
Clearly $\tau^{*}\equiv0$ if and only if $\sup_{\tau\in\mathcal
{T}}J(\tau)=u(s)$, which is equivalent to
$u(s)=\sup_{x\in(0,1]}[ w(x)u(\frac{s}{x})]$ by virtue of Theorem
\ref{uconvex}.

If $u(s)<\sup_{x\in(0,1]}[ w(x)u(\frac{s}{x})]$, then the last
inequality in \eqref{g2<g1} is strict unless $a=0$. This implies that
the maximum in \eqref{abstar} is not achievable.
\end{pf}

In the context of asset selling, $u(\cdot)$ is convex when the
underlying asset ranges from ``intermediate'' to ``bad'' depending on
the form of the original payoff
function $U(\cdot)$ (see the examples in Section~\ref{sec2}).
Theorem~\ref{uconvex} shows that in this case an optimal strategy is
\textit{in general} still of a ``take-profit-or-cut-loss'' form.
However, one must note that it is also possible that the maximum in
\eqref{abstar} is not achievable (as indicated in Corollary \ref
{care}). In that case
suppose $x^{*}=\argmax_{x\in(0,1]}[ w(x)u(\frac{s}{x})]$ exists. Let
$b^{*}=s/x^{*}$, and
\[
\tau_{(0, b^{*})}:=\inf\{t\geqslant 0\dvtx S_{t}\notin(0, b^{*})\}.
\]
Then we have
\[
\sup_{\tau\in\mathcal{T}}J(\tau)=\sup_{x\in(0,1]}\biggl[ w(x)u\biggl(\frac
{s}{x}\biggr)\biggr]=J\bigl(\tau_{(0, b^{*})}\bigr).
\]
However, when $b^{*}>s$, $\tau_{(0, b^{*})}$ is not a finite stopping
time [i.e., $P(\tau_{(0, b^{*})}=+\infty)>0$]. The interpretation of
this fact is that only a stop-gain threshold $b^*$ is set if one
applies $\tau_{(0, b^{*})}$; but as $S_t$ will never reach 0, with a
positive probability the process never exits the interval $(0,b^*)$.

\section{\texorpdfstring{Concave $u(\cdot)$}{Concave $u(.)$}}\label{sec5}
In this section, we study the case when $u(\cdot)$ is concave. Again,
we employ the quantile formulation \eqref{qua-form} where the
objective functional $J_{Q}(\cdot)$ becomes concave. In sharp contrast
to the case when $u(\cdot)$ is convex, in general the maxima of \eqref
{qua-form} are now in the interior of the constraint set, which can be
obtained using the classical Lagrange method. Let us first describe the
general solution procedure.

Consider a family of relaxed problems
%
\begin{eqnarray}\label{lagr}
J^{\lambda}_{Q}(G)&:=& \int_{0}^{1}u(G(x))w'(1-x)\,\dd x-\lambda\biggl(\int
_{0}^{1}G(x)\,\dd x-s\biggr)
\nonumber
\\
&=& \int_{0}^{1}\bigl(u(G(x))w'(1-x)-\lambda G(x)\bigr)\,\dd x+\lambda s\\
&=&\int
_{0}^{1}f^{\lambda}(x,G(x))\,\dd x +\lambda s,\nonumber
\end{eqnarray}
where $\lambda\geqslant 0$ and
$f^\lambda(x,y):=u(y)w'(1-x)-\lambda y.$
To maximize $J^{\lambda}_{Q}(\cdot)$ it suffices to maximize
$f^{\lambda}(x,\cdot)$ for each $x$. Recall that we do not assume
$u(\cdot)$ to be smooth [which is the case when, e.g., the original
payoff function $U(\cdot)$ is that of a call option]. Define
\begin{eqnarray*}
u'(x)&:=&\limsup_{h\to0+}\frac{u(x+h)-u(x)}{h},\\
(u')_{u}^{-1}(x)&:=&\inf\{y\geqslant 0\dvtx u'(y)<x\},\\
(u')_{l}^{-1}(x)&:=&\inf\{y\geqslant 0\dvtx u'(y)\leqslant x\}.
\end{eqnarray*}
It is easy to see that both $(u')_{l}^{-1}$ and $u'$ are
right-continuous, while $(u')_{u}^{-1}$ is left-continuous.
Fix $x$. As $f^{\lambda}(x,\cdot)$ is concave on $\R^+$, $y$
maximizes $f^{\lambda}(x,\cdot)$ on $\R^+$ if and only if
\[
y\in\biggl[(u')_{l}^{-1}\biggl(\frac{\lambda}{w'(1-x)}\biggr),(u')_{u}^{-1}\biggl(\frac
{\lambda}{w'(1-x)}\biggr)\biggr].
\]

To proceed, we need to further specify the shape of the probability
distortion function $w(\cdot)$. The case when $w(\cdot)$ is convex
has been solved in Section~\ref{seconvex} where $\tau^*= 0$ is an
optimal stopping time. The other cases will be studied in the next
three subsections, respectively.

\subsection{\texorpdfstring{Concave $w(\cdot)$}{Concave $w(.)$}}\label{sec5.1}

\begin{thmm}\label{uconcavewconcave}
If both $u(\cdot)$ and $w(\cdot)$ are concave, and there exists
$\lambda^{*}\geqslant 0$ such that $(u')_{l}^{-1}(\frac{\lambda
^{*}}{w'(1-x)})>0$, $\forall x\in(0,1)$ and
%
\begin{equation}\label{ld}
\int_{0}^{1}(u')_{l}^{-1}\biggl(\frac{\lambda^{*}}{w'(1-x)}\biggr)\,\dd x=s,
\end{equation}
then $G^*(x):=(u')_{l}^{-1}(\frac{\lambda^{*}}{w'(1-x)})$ is an
optimal solution to problem \eqref{qua-form}.
\end{thmm}

\begin{pf} Clearly $G^*(x)$
maximizes $f^{\lambda^*}({x},\cdot)$ on $\R^{+}$, for each $x\in(0,1)$.
Since $w'(1-x)$ is nondecreasing in $x$, $G^*$ is nondecreasing and
left-continuous. By defining $G^*(0)=0$ we see $G^*$ is indeed a
quantile function given that $G^*(x)>0$ $\forall x\in(0,1)$.
Moreover, $G^*\in\mathcal{Q}$ by virtue of \eqref{ld}.
On the other hand, for any \mbox{$G\in\mathcal{Q} $,}
\begin{eqnarray*}
J_{Q}(G)&\leqslant& J^{\lambda^*}_{Q}(G)=\int_{0}^{1}f^{\lambda
^*}(x,G(x))\,\dd x +\lambda^* s\\
&\leqslant&\int_{0}^{1}f^{\lambda^*}(x,G^*(x))\,\dd x +\lambda^*
s=J^{\lambda^*}_{Q}(G^*)=J_{Q}(G^{*}).
\end{eqnarray*}
So $G^*$ is an optimal solution to problem \eqref{qua-form}.
\end{pf}

The above general result involves an assumption that $\lambda^*$
exists so that \eqref{ld} holds. Conceptually, nonexistence of the Lagrange
multiplier is an indication of the ill-posedness or the
nonattainability of the underlying optimization problem (see Section 3
of~\cite{JXZ}
for a detailed study in the context of utility maximization).
Mathematically, when $u(\cdot)$ and $w(\cdot)$ are given
in specific forms (see, e.g., Example~\ref{ex1} below) it is straightforward to
check the validity of the assumption. In more general cases, one
constructs the function $\varphi(\lambda):=\int
_{0}^{1}(u')_{l}^{-1}(\frac{\lambda}{w'(1-x)})\,\mathrm{d}x$,\vspace*{1pt} and then checks
the validity of \eqref{ld} by examining the continuity of
$\varphi$ and its values at $\lambda=0$ and $\lambda\uparrow\infty$.

In general, the quantile function, $G^*$, of the optimally stopped
state does no longer correspond to a two-point distribution; or there
is no threshold level that would \textit{directly} trigger a stopping.
We have discussed in the previous section that $u(\cdot)$ being
concave corresponds to, at least in some cases of interest, an
``unfavorable'' underlying stochastic process. On the other hand, a
concave $w(\cdot)$ suggests that the agent is risk-seeking in that she
exaggerates the probability of the underlying process reaching a very
high state. In the context of stock selling, our result indicates that
a speculative agent, when holding
a ``bad'' stock, will not set any specific cut-loss or stop-gain prices.

Moreover, since $w(\cdot)$ is concave, we have
\[
\underline{b}:=(u')^{-1}_{l}\biggl(\frac{\lambda^*}{w'(1-)}\biggr)\leqslant
G_{\lambda^*}(x)\leqslant(u')^{-1}_{l}\biggl(\frac{\lambda
^*}{w'(0+)}\biggr)=:\bar b \qquad\forall x\in(0,1).
\]
In particular,
if $w'(0+)<\infty$, then $\bar b<+\infty$;
hence, the optimally stopped state will never exceed $\bar b$, or one
will have already stopped before the process ever reaches $\bar b$.
Similarly, if $w'(1-)>0$, then the optimally stopped state will never
fall below $\underline{b}$. If the range of $w'$ is a singleton which
must be $\{1\}$, then the range of the possible stopped states is also
a singleton, which is necessarily $\{s\}$. This shows that, in the case
of stock liquidation, if there is no probability distortion then
a bad stock will be sold immediately, which is also consistent with
Corollary~\ref{u-concave,w-convex}. In other words, if an agent is
still holding an unfavorable stock then it is an indication that the agent
is distorting probability scale hoping for extraordinarily return.

We now provide the following example to illustrate the general result
of Theorem~\ref{uconcavewconcave}.

\begin{exam}\label{ex1}
Consider a model of asset selling with a concave function $u(x)=\frac
{1}{\gamma}x^{\gamma}$, $0<\gamma<1$, and a concave distortion
function $w(x)=x^{\alpha}$, $0<\alpha\leqslant1$.
We have $(u')^{-1}(x)=x^{{1}/{(\gamma-1)}}$, $w'(1-x)=\alpha
(1-x)^{\alpha-1}$.

First we assume that $1>\alpha>\gamma$, namely, the agent is only
moderately risk-seeking (relative to the original payoff function and
the quality of the asset). The equation \eqref{ld} for $\lambda^*$ is
$\frac{1-\gamma}{\alpha-\gamma}(\frac{\lambda^*}{\alpha})^{
{1}/{(\gamma-1)}}=s$, which clearly has a unique solution.
Then the optimal quantile function is
%
\begin{equation}
G^{*}(x)=s\frac{\alpha-\gamma}{1-\gamma}\biggl(\frac{1}{1-x}\biggr)^{
{(1-\alpha)}/{(1-\gamma)}},\qquad x\in(0,1).
\end{equation}
The corresponding CDF of the optimally stopped price is
%
\begin{equation}
\qquad F^{*}(x)=\cases{
\displaystyle 1-\biggl(s\frac{\alpha-\gamma}{1-\gamma}\biggr)^{{(1-\gamma)}/{(1-\alpha
)}}x^{-{(1-\gamma)}/{(1-\alpha)}},& \quad $\displaystyle x\geqslant
s\frac{\alpha-\gamma}{1-\gamma};$\vspace*{2pt}\cr
\displaystyle 0,&\quad $\displaystyle x<s\frac{\alpha-\gamma}{1-\gamma}.$}
\end{equation}
This is a Pareto distribution\footnote{Pareto distribution was put
forth by Italian economist Vilfredo Pareto~\cite{Parato}
to describe the allocation of wealth among individuals in a society.}
with the Pareto index $\frac{1-\gamma}{1-\alpha}>1$.
In particular, one should never stop when the asset price is below
$s\frac{\alpha-\gamma}{1-\gamma}$, a true fraction of the initial
price $s$. Pareto index is a measure of the ``fatness'' of the tail of
the stopped price. The larger the Pareto index (i.e., the lower $\gamma
$ or the higher $\alpha$), the lighter tailed the distribution (and
hence, the smaller the proportion of very high stopped prices).
This makes perfect sense since a higher $\alpha$ implies a less
exaggeration of the probability of the asset achieving very high
prices, hence, more likely the agent stops at a moderate price.

There are infinitely many stopping times generating the same
distribution $F^{*}$ in this case. However, a convenient one is the
so-called\vadjust{\goodbreak} Az\'ema--Yor stopping time (see~\cite{AY})
%
\begin{equation}
\tau_{\mathrm{AY}}=\inf\biggl\{t\geqslant 0 \dvtx  S_t\leqslant\frac{\alpha-\gamma
}{1-\gamma}\max_{0\leqslant s\leqslant t}S_{s}\biggr\},
\end{equation}
which is an optimal solution to problem \eqref{objective}. Az\'
ema--Yor theorem is applicable in our case since
$\int_0^\infty x\,\dd F^*(x)\equiv\int_0^1G^*(x)\,\dd x=s$. Such a
stopping strategy is to stop at the first time when the boundary of the
\textit{drawdown constraint} $S_t\geqslant \frac{\alpha-\gamma
}{1-\gamma}\max_{0\leqslant s\leqslant t}S_{s}$
is touched upon. This implies that one sells as soon as the current
stock price falls below a true fraction of
the historical high price.

If $\alpha=1$ (i.e., there is no distortion), then $\tau_{\mathrm{AY}}=0$.
Hence, an agent who distorts the probability scale will hold an asset
which would be
otherwise sold immediately by one who does not. This shows that
probability distortion does change the optimal stopping behavior.

If $\alpha< \gamma$ so that the agent is sufficiently risk-taking,
then choose \textit{any} $0<\eta<1$ satisfying $\alpha<(1-\eta
)\gamma$. Take
$G^{*}(x)=\eta s(1-x)^{\eta-1}.$
It is easy to check that $G^{*}\in\mathcal{Q}$ while
\[
J_{Q}(G^{*})=\int_{0}^{1} \frac{1}{\gamma}\bigl(\eta s(1-x)^{\eta
-1}\bigr)^{\gamma} \alpha(1-x)^{\alpha-1}\,\dd x=+\infty.
\]
So the optimal value of \eqref{objective} in this case is $+\infty$
and $G^{*}$ is an optimal solution. Since $G^{*}$ also follows a Parato
distribution, the corresponding Az\'ema--Yor stopping time is given by
\[
\tau_{\mathrm{AY}}=\inf\Bigl\{t\geqslant 0 \dvtx  S_t\leqslant\eta\max_{0\leqslant
s\leqslant t}S_{s}\Bigr\}.
\]

Finally, when $\alpha= \gamma$ we construct
$G_{n}(x)=\frac{1}{n} s(1-x)^{{1}/{n}-1}, n>0.$
Then $G_{n}\in\mathcal{Q}$ with
\[
J_{Q}(G_{n})=\int_{0}^{1} \frac{1}{\gamma}\biggl(\frac{1}{n}
s(1-x)^{{1}/{n}-1}\biggr)^{\gamma} \alpha(1-x)^{\alpha-1}\,\dd x=\frac
{1}{\gamma}s^{\gamma}n^{1-\gamma}.
\]
Hence, the optimal value of the stopping problem is $+\infty$. The
corresponding Az\'ema--Yor stopping time is
\[
\tau_{{\mathrm{AY}},n}=\inf\biggl\{t\geqslant 0 \dvtx  S_t\leqslant\frac{1}{n}\max
_{0\leqslant s\leqslant t}S_{s}\biggr\}.
\]
It is not hard to show that there is no optimal solution in this case.
\end{exam}

\subsection{\texorpdfstring{Reverse $S$-shaped $w(\cdot)$}{Reverse $S$-shaped $w(.)$}}\label{subseu-concave,anti-s-dis}

\begin{thmm}\label{uconcavews}
Assume that $u(\cdot)$ is concave and $w(\cdot)$ is reverse
$S$-shaped, that is, it is concave on $[0,1-q]$ and convex on $[1-q,1]$
for some $q\in(0,1)$.
If $(a^*,\lambda^*)$ with $a^*>0$ is a solution to the following
mathematical program
%
\begin{eqnarray}\label{para}
\qquad&&\mbox{Maximize}\quad  \bigl(1-w(1-q)\bigr)u(a)\nonumber\\
&&\phantom{\mbox{Maximize}}\quad\qquad{}+\int_{q}^{1} u\biggl( a\vee
(u')_{l}^{-1}\biggl(\frac{\lambda}{w'(1-x)}\biggr)\biggr)w'(1-x)\,\dd x
\\
&&\quad\mbox{subject to:}\quad \lambda\geqslant 0,
a\geqslant0,
aq+\int_{q}^{1}a\vee(u')_{l}^{-1}\biggl(\frac{\lambda}{w'(1-x)}\biggr)\,\dd x =
s,\nonumber
\end{eqnarray}
then
%
\begin{eqnarray}\label{sshape}
G^{*}(x)=a^*\mathbf{1}_{(0,q]}(x)+\biggl(a^*\vee(u')^{-1}_{l}\biggl(\frac
{\lambda^*}{w'(1-x)}\biggr)\biggr)\mathbf{1}_{(q,1)}(x),
\nonumber
\\[-8pt]
\\[-8pt]
\eqntext{  x\in[0,1)}
\end{eqnarray}
is an optimal solution to problem \eqref{qua-form}.
\end{thmm}

The proof of this theorem is rather technical, and it is delayed to
Appendix~\ref{AE}.

Reverse $S$-shaped probability distortion has been used and studied by
many authors (see, e.g.,~\cite{TverskyAFoxC95pd,PrelecD98df,Jin&Zhou})
and in particular by Kahneman and Tversky in the celebrated
CPT~\cite{T&K1992}. For a reverse $S$-shaped distortion $w(\cdot)$,
$w'(x)>1$ around both $x=0$ and $x=1$. This implies, as seen from
\eqref{qua-form}, that such distortion puts higher weights on \textit
{both} very good and very bad outcomes. In other words, the agent
exaggerates the small probabilities of both very good and very bad scenarios.
In~\cite{HZ}, exaggeration of small probabilities for extremely good
and bad outcomes is used to model the emotion of \textit{hope} and
\textit{fear}, respectively. In the current context of optimal
stopping, the expression \eqref{sshape} shows, qualitatively, that the
agent sets a cut-loss
level $a$ (because she has fear) and does not set any stop-gain level
(because she has hope). This is widely known as the
``cut-loss-and-let-profit-run''
strategy in stock trading.

\begin{exam}\label{ex2}
Consider a concave $u(x)=\frac{1}{\gamma}x^{\gamma}$, $0<\gamma<1$,
and a reverse $S$-shaped distortion function
\[
w(x)=\cases{
 2x-2x^{2}, & \quad$ 0\leqslant x\leqslant\frac{1}{2};$\vspace*{2pt}\cr
2x^{2}-2x+1, &\quad $\frac{1}{2}< x\leqslant1.$}
\]
Then constraints in \eqref{para} become
%
\[
\lambda\geqslant 0, a\geqslant0,\qquad
\frac{1}{2}a+\int_{{1}/{2}}^{1}a\vee\biggl(\frac{\lambda
}{4x-2}\biggr)^{{1}/{(\gamma-1)}}\,\dd x= s
\]
and the objective function in \eqref{para} is
%
\[
J(a,\lambda)= \frac{1}{2\gamma}a^{\gamma}+\frac{1}{\gamma}\int
_{{1}/{2}}^{1}a^{\gamma}\vee\biggl(\frac{\lambda}{4x-2}\biggr)^{
{\gamma}/{(\gamma-1)}}(4x-2)\,\dd x.\vadjust{\goodbreak}
\]
Define
\[
\bar{c}=\inf\biggl\{x\geqslant \frac{1}{2}\dvtx a\leqslant\biggl(\frac{\lambda
}{4x-2}\biggr)^{{1}/{(\gamma-1)}}\biggr\}\wedge1\in[0.5,1].
\]

If $\bar{c}=1$, then the problem reduces to maximize $J(a,\lambda)=
\frac{1}{\gamma}s^{\gamma}$ subject to
$\lambda\geqslant 0$, $a= s$, which is trivial.
If $\bar{c}\in(0.5,1)$, then the above constraints are equivalent to
%
\begin{eqnarray}\label{oppara}
a&=&\frac{s}{\bar{c}+{(1-\gamma)}/({2(2-\gamma)}) ((2\bar
{c}-1)^{{1}/{(\gamma-1)}}-(2\bar{c}-1) )},
\nonumber
\\[-8pt]
\\[-8pt]
\nonumber
\lambda&=&(4\bar{c}-2)a^{\gamma-1},\qquad
0.5< \bar{c}<1,
\end{eqnarray}
and thus our objective is to maximize
$
\frac{1}{\gamma}s^{\gamma} g(\bar{c})
$
in $\bar{c}\in(0.5,1)$, where
\begin{eqnarray*}
g(\bar{c})&:=&\biggl(\frac{1}{\bar{c}+{(1-\gamma)}/({2(2-\gamma)})
((2\bar{c}-1)^{{1}/{(\gamma-1)}}-(2\bar{c}-1) )}\biggr)^{\gamma}
\\
&&{}\times\biggl(1-2\bar{c}+2\bar{c}^{2}+\frac{1-\gamma}{2-\gamma}\bigl((2\bar
{c}-1)^{{\gamma}/{(\gamma-1)}} -(2\bar{c}-1)^{2}\bigr)\biggr).
\end{eqnarray*}

Now, $g(0.5+)=(\frac{1-\gamma}{2-\gamma})^{1-\gamma}2^{\gamma
}=(\frac{1-\gamma}{4(2-\gamma)}2^{{(2-\gamma)}/{(1-\gamma
)}})^{1-\gamma}$ and $g(1-)=g(1)=1$.
Noting
$2^{t}<4t\ \forall t\in[2,4)$,
we conclude
$\frac{1-\gamma}{4(2-\gamma)}2^{{(2-\gamma)}/{(1-\gamma)}}<1$,
if $2<\frac{2-\gamma}{1-\gamma}<4$ or $0<\gamma<\frac{2}{3}.$
Thus
$g(0.5+)<g(1-)$, if $0<\gamma<\frac{2}{3}.$
If $\bar{c}=0.5$, then $a=0$ and the optimal value is  $\frac{1}{\gamma}s^{\gamma} g(0.5+)$.
In other words,
the maximum value of the objective function is achieved at some point $\bar{c}^*\in (0.5,1]$.

We have now deduced the optimal quantile function
\begin{eqnarray*}
G^{*}(x)=a\mathbf{1}_{(0,\bar{c}]}(x)+ a \biggl(\frac{4\bar
{c}-2}{4x-2}\biggr)^{{1}/{(\gamma-1)}}\mathbf{1}_{(\bar{c},1)}(x),\qquad x\in[0,1),
\end{eqnarray*}
where $\bar{c}\equiv\bar{c}^*\in(0.5,1]$, and $a>0$ is determined
via \eqref{oppara}. The corresponding optimal CDF is
\begin{eqnarray*}
F^{*}(x)=\cases{
0,&\quad $x<a;$\vspace*{2pt}\cr
\displaystyle (\bar{c}-0.5)(x/a)^{1-\gamma}+0.5, &\quad $a\leq x<(2\bar
{c}-1)^{{1}/{(\gamma-1)}}a
;$\vspace*{2pt}\cr
1,&\quad $ x\geqslant(2\bar{c}-1)^{{1}/{(\gamma-1)}}a.$}
\end{eqnarray*}
The barycenter function (also called the Hardy--Littlewood maximal function) for a \textit{centered} probability
measure $F$  is generally given by
\[
\Psi_F(x)=\cases{
0,&\quad  $x\leq m;$ \vspace*{2pt}\cr
\displaystyle\frac{\int_{[x,\infty)}y\,\dd F(y)}{1-F(x-)}, &\quad $ m< x<M
;$\vspace*{2pt}\cr
 x,&\quad  $x\geqslant M,$}
\]
where
$m=\inf\{x \dvtx F(x)>0\}$ and $M=\sup\{x \dvtx F(x)<1\}$ (see  Az\'ema and Yor~\cite{AY}).
In our case, it is
\[
\Psi(x)
=\cases{
s,\qquad x\leq a; \vspace*{2pt}\cr
\displaystyle\frac{1-\gamma}{2-\gamma}
 \frac{(x/a)^{2-\gamma}-(2\bar{c}-1)^{{(2-\gamma)}/{(1-\gamma)}}}{(x/a)^{1-\gamma}-(2\bar{c}-1)}a,\vspace*{2pt}\cr
\hspace*{33pt}a< x<(2\bar{c}-1)^{{1}/{(\gamma-1)}}a ;\vspace*{2pt}\cr
x,\qquad x\geq (2\bar{c}-1)^{{1}/{(\gamma-1)}}a,}
\]
whereas the corresponding Az\'ema--Yor stopping time is
\[
\tau_{\mathrm{AY}}=\inf\Bigl\{t\geqslant 0 \dvtx  \Psi(S_t)\leqslant \max
_{0\leqslant s\leqslant t}S_{s}\Bigr\}.
\]

Suppose now that $\gamma=0.3$. Then the optimal $\bar{c}^* \approx
0.70$, and $a \approx0.72s$, $\lambda\approx s^{-0.7}$ are determined
via \eqref{oppara}. Therefore the optimal quantile function presented
in Theorem~\ref{uconcavews} is
\[
G^{*}(x) \approx0.72s\mathbf{1}_{(0,0.7]}(x)+ (4x-2)^{1.43}s\mathbf
{1}_{(0.7,1)}(x).
\]

Since the barycenter function $\Psi$ is increasing, the Az\'ema--Yor
stopping time is the first time when $S_t$ hits $\Psi^{-1}(\max
_{0\leqslant s\leqslant t}S_{s})$, a moving level that is related to
the running maximum (rather than a proportion of the running maximum as
in Example~\ref{ex1}).
\end{exam}

\subsection{\texorpdfstring{$S$-shaped $w(\cdot)$}{$S$-shaped $w(.)$}}\label{subseu-concave,s-dis}

\begin{thmm}\label{uconcavewsr}
Assume that $u(\cdot)$ is concave, and $w(\cdot)$ is $S$-shaped, that
is, it is convex on $[0,1-q]$ and concave on $[1-q,1]$ for some $q\in(0,1)$.
If $(a^*,\lambda^*)$ is a solution to the following mathematical program
\begin{eqnarray*}
&&\mbox{Maximize}\quad  \int_{0}^{q} u\biggl(a\wedge(u')_{l}^{-1}\biggl(\frac
{\lambda}{w'(1-x)}\biggr)\biggr)w'(1-x)\,\dd x +w(1-q) u(a)\\
&&\quad\mbox{subject to:}\quad \lambda\geqslant 0,
a\geqslant0,
\int_{0}^{q}a\wedge(u')_{l}^{-1}\biggl(\frac{\lambda}{w'(1-x)}\biggr)\,\dd x
+a(1-q) = s
\end{eqnarray*}
and $(u')_{l}^{-1}(\frac{\lambda^*}{w'(1-x)})>0$ $\forall x\in(0,q]$,
then
\begin{eqnarray*}
G^{*}(x)=\biggl(a^*\wedge(u')_{l}^{-1}\biggl(\frac{\lambda^*}{w'(1-x)}\biggr)\biggr)\mathbf
{1}_{(0,q]}(x)+a^*\mathbf{1}_{(q,1)}(x),\qquad x\in[0,1),
\end{eqnarray*}
is an optimal solution to problem \eqref{qua-form}.
\end{thmm}

\begin{pf} The proof is similar to that of Theorem~\ref{uconcavews};
hence, omitted.~%
\end{pf}

The economic interpretation of the result for this case is just
opposite to the reverse $S$-shaped counterpart. An $S$-shaped
probability distortion
underlines an agent\vadjust{\goodbreak} who under-weighs probabilities of extreme events
(both good and bad). So she sets an upper target level simply because
she is not hopeful
for a dramatically high price, while she does not prescribe a cut-loss
level since she believes the asset will not go catastrophically wrong.

\subsection{Discussion}\label{sec5.4}

We have obtained the quantile functions of the optimally stopped states
for the three cases discussed in this section. In order to finally
solve the original distorted optimal stopping problem \eqref
{objective}, we need to recover optimal stopping times from the
quantile functions. Unlike the cases investigated in Section \ref
{seconvex} where
the optimal distribution/quantile functions are those of two-point or
one-point distributions and optimal stopping times can be uniquely
determined, there could be infinitely many stopping times corresponding
to the same distribution of the stopped state.\footnote{In the
Skorokhod embedding literature,
one usually introduces additional criteria in order to uniquely
determine the stopping time (see, e.g.,~\cite{Jan}).} As demonstrated
in Examples~\ref{ex1} and~\ref{ex2}, the Az\'ema--Yor stopping time would provide
a convenient solution that is related to the running maximum of the
underlying process, which is also commonly incorporated in practice. On
the other hand, for many applications, how optimally stopped states are
probabilistically distributed already reveals important qualitative
information. For instance, we have shown in this section when the agent
would put a cut-loss floor or a target state or simply set none,
depending on her risk preferences. An optimally stopped state
distribution is also adequate in calculating the optimal payoff value
function, which is
relevant in the context of, say, option pricing or irreversible investment.

The results in this section have also demonstrated how probability
distortion affects optimal stopping strategies. In the previous section
we have proved that if there is no probability distortion,
optimal stopping strategies are always of the threshold-type with at
most two thresholds. Thus one stops only at (at most) two states.
Strategies are qualitatively changed when there is probability
distortion, where one sets only one-sided threshold or simply none.
Moreover, there could be infinitely many stopped states.

\section{\texorpdfstring{$S$-shaped $u(\cdot)$}{$S$-shaped $u(.)$}}\label{sec6}
We now consider the case when $u(\cdot)$ is $S$-shaped. If
the distortion $w(\cdot)$ is convex, then the result has already been
derived in Section~\ref{seconvex}. If $w(\cdot)$ is concave, then
we can utilize the same idea as in Section~\ref{subseu-concave,s-dis} to get similar results,
thanks to the duality between $u(\cdot)$ and $w(\cdot)$. If $w(\cdot
)$ is $S$-shaped, we can also apply similar techniques. We leave the
details to the interested readers. In this section, we will focus on
the most interesting case when $w(\cdot)$ is a reverse $S$-shaped
distortion function.\footnote{If $u(\cdot)$ is interpreted as a
utility function, then
the case when $u(\cdot)$ is $S$-shaped while $w(\cdot)$ is reverse
$S$-shaped is consistent with the CPT
of~\cite{T&K1992}.}

Henceforth in this section $u$ is convex on $[0,\theta]$ and concave
on $[\theta,\infty)$ for some $\theta>0$, and
$w$ is concave on $[0,1-q]$ and convex on $[1-q,1]$ for some $q\in(0,1)$.

Fix $G_{0}\in\mathcal{G}$.
Let
$x_{0}=\sup\{x\in[0,1) | G_{0}(x)\leqslant\theta\}\wedge1.$
Then $G_{0}(x_{0})\leqslant\theta$ since $G_{0}$ is left-continuous, and
\begin{eqnarray*}
J_{Q}(G_{0})=\int_{0}^{x_{0}}u(G_{0}(x))w'(1-x)\,\dd x +\int
_{x_{0}}^{1}u(G_{0}(x))w'(1-x)\,\dd x.
\end{eqnarray*}

Consider two subproblems:
%
\begin{eqnarray}
&\label{g-}\displaystyle\max_{G\in\mathcal{G}^{-} }\int_{0}^{x_{0}}u(G(x))
w'(1-x)\,\dd x,&\\
&\label{g+}\displaystyle\max_{G\in\mathcal{G}^{+} }\int_{x_{0}}^{1}u(G(x))
w'(1-x)\,\dd x,&
\end{eqnarray}
where
\begin{eqnarray*}
\mathcal{G}^{-}&=&\biggl\{G\in\mathcal{G} \Big| G(x_{0})\leqslant
G_{0}(x_{0}),\int_{0}^{x_{0}}G(x)\,\dd x\leqslant\int
_{0}^{x_{0}}G_{0}(x) \,\dd x\biggr\},\\
\mathcal{G}^{+}&=&\biggl\{G\in\mathcal{G} \Big| G(x_{0}+)\geqslant
G_{0}(x_{0}) ,\int_{x_{0}}^{1}G(x)\,\dd x\leqslant\int
_{x_{0}}^{1}G_{0}(x) \,\dd x\biggr\}.
\end{eqnarray*}
Subproblem \eqref{g-} is a convex maximization problem. Using the idea
of the proof of Lemma~\ref{2by1}, we can show that the optimal
solution to the subproblem \eqref{g-} is of the form
\begin{eqnarray*}
G(x)=a_1\mathbf{1}_{(0,c_{1}]}(x)+a_2\mathbf
{1}_{(c_{1},c_{2}]}(x)+G_{0}(x_{0})\mathbf{1}_{(c_{2},x_{0}]}(x)\qquad
\forall x\in(0,x_{0}].
\end{eqnarray*}
For subproblem \eqref{g+}, we can use the idea of proof of Theorem
\ref{uconcavews} to show that the optimal solution must be of the form
\begin{eqnarray}
G(x)=G_{0}(x_{0})\mathbf{1}_{(x_{0},q]}(x)+\biggl(G_{0}(x_{0})\vee
(u')^{-1}_{l}\biggl(\frac{\lambda}{w'(1-x)}\biggr)\biggr)\mathbf{1}_{(q,1)}(x)
\nonumber
\\
 \eqntext{\forall x\in(x_{0},1).}
\end{eqnarray}
Now we conclude that the optimal solution is of the form
\begin{eqnarray*}
G^{*}(x)&=&a_1\mathbf{1}_{(0,c_{1}]}(x)+a_2\mathbf
{1}_{(c_{1},c_{2}]}(x)+a_3\mathbf{1}_{(c_{2},q]}(x) \\
&&{}+\biggl(a_3\vee
(u')^{-1}_{l}\biggl(\frac{\lambda}{w'(1-x)}\biggr)\biggr)\mathbf{1}_{(q,1)}(x)\qquad
\forall x\in(0,1),
\end{eqnarray*}
where parameters $a_1$, $a_2$, $a_3$, $c_1$, $c_2$ and $\lambda$ are
subject to
%
\begin{eqnarray*}
&&\lambda\geqslant 0,
0< a_1\leqslant a_2\leqslant a_3\leqslant\theta ,
0\leqslant c_{1}\leqslant c_{2}\leqslant q,\\
&&\qquad a_1c_1+a_2(c_2-c_1)+a_3(q-c_2)+\int_{q}^{1} a_3\vee
(u')^{-1}_{l}\biggl(\frac{\lambda}{w'(1-x)}\biggr)\,\dd x\leqslant s.
\end{eqnarray*}
The objective is
\begin{eqnarray*}
J_{Q}(G^{*})&=&\bigl(1-w(1-c_1)\bigr)u(a_1)+\bigl(w(1-c_1)-w(1-c_2)\bigr)u(a_2)\\
&&{}+\bigl(w(1-c_2)-w(1-q)\bigr)u(a_3)\\
&&{}+\int_{q}^{1} u\biggl(a_3\vee(u')_{l}^{-1}\biggl(\frac{\lambda
}{w'(1-x)}\biggr)\biggr)w'(1-x)\,\dd x.
\end{eqnarray*}
Hence, the original problem reduces to the above mathematical program
which can be solved readily.

\section{Concluding remarks}\label{sec7}

In this paper we have formulated an optimal stopping problem under
distorted probabilities and developed an approach, primarily based on
the distribution/quantile formulation and the Skorokhod embedding, to
solving this new problem. Note that the optimal stopping strategies
derived are \textit{pre-committed}, instead of dynamically consistent.
Precisely, while our solutions are optimal at $t=0$, they are no longer
optimal at $t=\varepsilon$ for any $\varepsilon>0$. This is due to
the inherent time-inconsistency arising from the distortion. There are
recent studies on time-inconsistent optimal control, which use a
time-consistent game equilibrium to replace the notion of
``optimality'' (see, e.g.,~\cite{EL} and
\cite{BMZ}). It is, however, not clear how to extend this equilibrium
idea to the optimal stopping setting.
On the other hand, a pre-committed strategy is still important. It will
determine the \textit{value} of the problem at any given time.
Moreover, in reality people often uphold strategies for a certain time
period before changing them, even though they are dealing with
time-inconsistent problems which may call for continuously changing
strategies. For example, Barberis~\cite{B} analyzed in detail, in the
setting of casino gambling, the behavior of a ``sophisticated''
gambler who can commit to his initial exit strategy.

An important point to note is that, while our stopping strategies are
time-inconsistent in terms of
quantitative values, they are indeed time-consistent in terms of
qualitative types. For example, Theorem~\ref{convexw} stipulates that if the current
$S_t=s$, then the optimal stopping strategy is to stop either at $a^*$
or $b^*$, whose values depend on $s$ via (\ref{abstar}). At the next moment the
underlying process value becomes $S_{t+\varepsilon}=s'$, then one
needs to re-calculate (\ref{abstar}) to obtain the new thresholds $a'$ and $b'$.
Although the strategy has changed, its type (that of two-threshold) has
not, which depends only on the risk preference of the agent and the
property of the process.

We assume in this paper the underlying stochastic process to be a GBM
for two reasons: (1) it is widely used in many applications especially
in finance and (2) we would like to concentrate on the new approach
developed (which is already very complex) without being carried away by
the complexity of a more general underlying process. The advantage of a
GBM is that it can be turned into an exponential martingale via a
simple transformation; thus the Skorokhod\vadjust{\goodbreak} embedding applies. A~more
general process governed by a nonlinear SDE may still be transformed
into a martingale, but more technicalities need to be taken care of,
especially in terms of the range of the martingale. This will be
studied in a forthcoming paper.

\begin{appendix}


\section{\texorpdfstring{Proof of Theorem \lowercase{\protect\ref{bah}}}{Proof of Theorem 2.1}}\label{AA}

Because $u(\cdot)$ is nonincreasing, we have
\begin{eqnarray*}
\sup_{\tau\in\mathcal{T}}J(\tau)=\sup_{\tau\in\mathcal{T}}\int
_{0}^{u(0+)}w\bigl(\BP\bigl(u(S_{\tau})> x\bigr)\bigr)\,\dd x\leqslant\sup_{\tau\in
\mathcal{T}}\int_{0}^{u(0+)}1\,\dd x=u(0+).
\end{eqnarray*}
On the other hand,
\begin{eqnarray*}
\sup_{\tau\in\mathcal{T}}J(\tau)&\geqslant& \limsup_{T\to+\infty
}J(T)\geqslant \liminf_{T\to+\infty}J(T)\geqslant
\liminf_{T\to+\infty}\int_{0}^{\infty}w\bigl(\BP\bigl(u(S_{T})> x\bigr)\bigr)\,\dd x\\
&\geqslant& \int_{0}^{\infty}\liminf_{T\to+\infty}w\bigl(\BP\bigl(u(S_{T})>
x\bigr)\bigr)\,\dd x
\geqslant \int_{0}^{\infty}w\Bigl(\liminf_{T\to+\infty}\BP\bigl(u(S_{T})>
x\bigr)\Bigr)\,\dd x\\
&=&\int_{0}^{\infty}w\bigl(\BP\bigl(u(0+)\geqslant x\bigr)\bigr)\,\dd x=u(0+)\geqslant
\sup_{\tau\in\mathcal{T}}J(\tau),
\end{eqnarray*}
where we used the fact that $\lim_{t\to\infty} S_{t}=0$ almost
surely since $S_t$ is an exponential martingale.
This implies that $u(0+)$ is the optimal value of problem \eqref
{objective}, and \eqref{infinite} holds.

Next, the fact that $\lim_{t\to\infty} S_{t}=0$ implies $\tau_{\ell
}\in\mathcal{T} $ for any $\ell>0$.
Now, if there is $\ell>0$ such that $u(\ell)=u(0+)$, then
$u(x)=u(0+)$ for each $x\in(0,\ell)$ since $u(\cdot)$ is
nonincreasing and therefore
$
\sup_{\tau\in\mathcal{T}}J(\tau)\leqslant u(0+)=J(\tau_{\ell}),
$
proving that $\tau_{\ell}$ solves problem \eqref{objective}.

If there is no $\ell>0$ such that $u(\ell)=u(0+)$, then for every
fixed $\tau\in\mathcal{T}$, we have $u(S_{\tau})<u(0+)$ almost
surely. Consequently,
\begin{eqnarray*}
J(\tau)=\int_{0}^{u(0+)}w\bigl(\BP\bigl(u(S_{\tau})> x\bigr)\bigr)\,\dd x<\int
_{0}^{u(0+)}w\bigl(\BP\bigl(u(0+)> x\bigr)\bigr)\,\dd x=u(0+),
\end{eqnarray*}
which shows that there is no optimal solution to problem \eqref{objective}.


\section{\texorpdfstring{Proof of Lemma \lowercase{\protect\ref{2forms}}}{Proof of Lemma 3.1}}\label{AB}

Let $F$ and $G$ be the CDF and the quantile function of $S_{\tau}$,
respectively, for a fixed $\tau\in\mathcal{T}$.

First we assume that $u(\cdot)$ is a strictly increasing, $C^\infty$
function with $u(0)=0$. Then
\begin{eqnarray*}
J(\tau)&=&\int_{0}^{\infty}w\bigl(\BP\bigl(u(S_{\tau})>x\bigr)\bigr)\,\dd x
=\int_{0}^{\infty}w\bigl(\BP\bigl(u(S_{\tau})>u(y)\bigr)\bigr)\,\dd u(y)\\
&=&\int_{0}^{\infty}w\bigl(\BP(S_{\tau}>x)\bigr)\,\dd u(x)
= \int_{0}^{\infty}w\bigl(1-F(x)\bigr)\,\dd u(x)\\
&=&\int_{0}^{\infty}u(x)\,\dd\bigl[- w\bigl(1-F(x)\bigr)\bigr]
=\int_0^{\infty}u(x)w'\bigl(1-F(x)\bigr)\,\dd F(x)\\
&=&\int_0^1u(G(x))w'(1-x)\,\dd x,
\end{eqnarray*}
which proves \eqref{qua-form}, where the fifth equality is due to
Fubini's theorem.

Next, given an absolutely continuous, nondecreasing function $u(\cdot
)$ with $u(0)=0$, for each $\varepsilon>0$, we can find a strictly
increasing, $C^\infty$ function $u_{\varepsilon}(\cdot)$ such that
$|u_{\varepsilon}(x)-u(x)|<\varepsilon$, for all $ x\in\R^{+}.$
It is easy to check that
\begin{eqnarray*}
\biggl|\int_{0}^{\infty}w\bigl(\BP\bigl(u_{\varepsilon}(S_{\tau})>x\bigr)\bigr)\,\dd x-\int
_{0}^{\infty}w\bigl(\BP\bigl(u(S_{\tau})>x\bigr)\bigr)\,\dd x\biggr|&\leqslant&\varepsilon,\\
\biggl|\int_0^1u_{\varepsilon}(G(x))w'(1-x)\,\dd x-\int_0^1u(G(x))w'(1-x)\,\dd
x\biggr|& \leqslant&\varepsilon.
\end{eqnarray*}
We have proved that
\[
\int_{0}^{\infty}w\bigl(\BP\bigl(u_{\varepsilon}(S_{\tau})>x\bigr)\bigr)\,\dd x=\int
_0^1u_{\varepsilon}(G(x))w'(1-x)\,\dd x.
\]
Therefore,
\[
\biggl|\int_0^1u(G(x))w'(1-x)\,\dd x-\int_{0}^{\infty}w\bigl(\BP\bigl(u(S_{\tau
})>x\bigr)\bigr)\,\dd x\biggr| \leqslant 2\varepsilon.
\]
Since $\varepsilon$ is arbitrary, \eqref{qua-form} follows.

To show \eqref{dis-form}, we note
\begin{eqnarray*}
J(\tau)&=&\int_0^1u(G(x))w'(1-x)\,\dd x=\int_0^1u(G(x))\,\dd[-w(1-x)]\\
&=&\int_{0}^{\infty}u(x)\,\dd\bigl [- w\bigl(1-F(x)\bigr)\bigr]=\int_0^{\infty
}w\bigl(1-F(x)\bigr)\,\dd u(x)
\\
&=&\int_{0}^{\infty}w\bigl(1-F(x)\bigr)u'(x)\,\dd x,
\end{eqnarray*}
where the fourth equality is due to Fubini's theorem.

Finally, \eqref{equiv} is evident.

\section{\texorpdfstring{Proof of Lemma \lowercase{\protect\ref{f=f2}}}{Proof of Lemma 4.1}}\label{AC}
To prove this lemma we need some technical preliminaries.

\begin{lemma}\label{2by1}
For any $F^{*}\in\mathcal{S}_{n+1}$, $n=2,3,\ldots,$ there exist
$F_{1}$, $F_{2}\in\mathcal{S}_{n}$ and $\theta\in[0,1]$ such that
%
\begin{eqnarray*}
F^{*}&=&\theta F_{1}+(1-\theta)F_{2},\\
\int_{0}^{\infty}\bigl(1-F_{1}(x)\bigr)\,\dd x&=&\int_{0}^{\infty}\bigl(1-F_{2}(x)\bigr)\,\dd
x=\int_{0}^{\infty}\bigl(1-F^{*}(x)\bigr)\,\dd x.
\end{eqnarray*}
\end{lemma}

\begin{pf}
We first prove the lemma for $n=2$. Suppose $F^{*}\in\mathcal
{S}_{3}$. Write
\[
F^{*}= c_{1}\mathbf{1}_{[a_{1}, a_{2})}+ c_{2}\mathbf{1}_{[ a_{2},
a_{3})}+\mathbf{1}_{[ a_{3},\infty)},\qquad
s_{0}:=\int_{0}^{\infty}\bigl(1-F^{*}(x)\bigr)\,\dd x,
\]
where $a_1<a_2<a_3$ (otherwise the desired result holds trivially).
Note that
\[
a_{1}=\int_{0}^{a_{1}}\bigl(1-F^{*}(x)\bigr)\,\dd x\leqslant\int_{0}^{\infty
}\bigl(1-F^{*}(x)\bigr)\,\dd x= \int_{0}^{ a_{3}}\bigl(1-F^{*}(x)\bigr)\,\dd x\leqslant a_{3};
\]
that is, $a_{1}\leqslant s_{0}\leqslant a_{3}$. If $s_{0}=a_{1}$, or
$s_{0}= a_{3}$, then $F^{*}\in\mathcal{S}_{2}$ and we are done by
choosing $F_{1}=F_{2}=F^{*}$. Hence, from now on, we assume
$a_{1}<s_{0}< a_{3}$.

If $s_{0}> a_{2}$, then let
$F_{1}:=b_{1}\mathbf{1}_{[a_{1}, a_{3})}+\mathbf{1}_{[ a_{3},\infty
)}$ and
$F_{2}:=b_{2}\mathbf{1}_{[ a_{2}, a_{3})}+\mathbf{1}_{[ a_{3},\infty)},$
where
$b_{1}=\frac{ a_{3}-s_{0}}{ a_{3}-a_{1}}\in(0,1)$ and $b_{2}=\frac{
a_{3}-s_{0}}{ a_{3}- a_{2}}\in(0,1).$
It follows from
%
\[
s_{0}\equiv\int_{0}^{\infty}\bigl(1-F^{*}(x)\bigr)\,\dd x=a_{1}+(
a_{2}-a_{1})(1- c_{1})+( a_{3}- a_{2})(1- c_{2})
\]
that
$\frac{ c_{1}}{b_{1}}+\frac{ c_{2}- c_{1}}{b_{2}}=1.$
It is now easy to see that $F_{1}$, $F_{2}$ and $\theta:={
c_{1}}/{b_{1}}$ satisfy the desired requirements.

If $s_{0}\leqslant a_{2}$, then let
$F_{1}:=b_{1}\mathbf{1}_{[a_{1}, a_{3})}+\mathbf{1}_{[ a_{3},\infty
)}$ and
$F_{2}:=b_{2}\mathbf{1}_{[a_{1}, a_{2})}+\mathbf{1}_{[ a_{2},\infty)},$
where
$b_{1}=\frac{ a_{3}-s_{0}}{ a_{3}-a_{1}}\in(0,1)$ and $b_{2}=\frac{
a_{2}-s_{0}}{ a_{2}-a_{1}}\in[0,1).$
Define
$\theta_{1}:=\frac{ c_{1}- c_{2}b_{2}}{b_{1}(1-b_{2})}$, and
$\theta_{2}:=\frac{ c_{2}- c_{1}}{1-b_{2}}\geqslant 0.$
Noting that
\begin{eqnarray*}
s_{0}&\equiv&\int_{0}^{\infty}\bigl(1-F^{*}(x)\bigr)\,\dd x=a_{1}+(
a_{2}-a_{1})(1- c_{1})+( a_{3}- a_{2})(1- c_{2})\\
&\geqslant& a_{1}+(
c_{2}- c_{1})( a_{2}-a_{1}),
\end{eqnarray*}
we deduce $\theta_{2}\leqslant1$.
It is an easy exercise to verify
\[
\theta_{1}(1-F_{1})+\theta_{2}(1-F_{2})-(\theta_{1}+\theta
_{2})\mathbf{1}_{[0, a_{3})}=1-F^{*}-\mathbf{1}_{[0, a_{3})}.
\]
Integrating both sides on $(0,\infty)$, we obtain
$\theta_{1}s_{0}+\theta_{2}s_{0}-(\theta_{1}+\theta_{2})
a_{3}=s_{0}- a_{3},$
which leads to $\theta_{1}+\theta_{2}=1$ noting $s_{0}< a_{3}$.
Now we can easily verify that $F_{1}$, $F_{2}$ and $\theta:=\theta
_{1}$ satisfy the desired properties.

Now, let $F^{*}\in\mathcal{S}_{n+1}$ where $n=3,4,\ldots.$ Write
\begin{eqnarray*}
F^{*}&=& c_{1}\mathbf{1}_{[a_{1}, a_{2})}+ c_{2}\mathbf{1}_{[ a_{2},
a_{3})}+ c_{3}\mathbf{1}_{[ a_{3},a_{4})}+\sum_{i=4}^{n-1}
c_{i}\mathbf{1}_{[a_{i},a_{i+1})}+\mathbf{1}_{[a_{n},\infty)}\\
&=& c_{3}\mathbf{1}_{[a_{1},a_{4})}\biggl(\frac{ c_{1}}{ c_{3}}\mathbf
{1}_{[a_{1}, a_{2})}+\frac{ c_{2}}{ c_{3}}\mathbf{1}_{[ a_{2},
a_{3})}+\mathbf{1}_{[ a_{3},\infty)}\biggr)+\sum_{i=4}^{n-1} c_{i}\mathbf
{1}_{[a_{i},a_{i+1})}+\mathbf{1}_{[a_{n},\infty)}.
\end{eqnarray*}
Denote
\[
\bar{F}=\frac{ c_{1}}{ c_{3}}\mathbf{1}_{[a_{1}, a_{2})}+\frac
{ c_{2}}{ c_{3}}\mathbf{1}_{[ a_{2}, a_{3})}+\mathbf{1}_{[
a_{3},\infty)}\in\mathcal{S}_{3}.
\]
By what we have proved above, there exist $\bar{F}_{1}$,
$\bar{F}_{2}\in\mathcal{S}_{2}$ and $\theta\in[0,1]$ such that
%
\begin{eqnarray*}
\bar{F}&=&\theta\bar{F}_{1}+(1-\theta)\bar{F}_{2},\\
\int_{0}^{\infty}\bigl(1-\bar{F}_{1}(x)\bigr)\,\dd x&=&\int_{0}^{\infty
}\bigl(1-\bar{F}_{2}(x)\bigr)\,\dd x=\int_{0}^{\infty}\bigl(1-\bar
{F}(x)\bigr)\,\dd x.
\end{eqnarray*}
Define
\begin{eqnarray*}
F_{1}&:=& c_{3}\mathbf{1}_{[a_{1},a_{4})}\bar{F}_{1}+\sum
_{i=4}^{n-1} c_{i}\mathbf{1}_{[a_{i},a_{i+1})}+\mathbf
{1}_{[a_{n},\infty)},\\
F_{2}&:=& c_{3}\mathbf{1}_{[a_{1},a_{4})}\bar{F}_{2}+\sum
_{i=4}^{n-1} c_{i}\mathbf{1}_{[a_{i},a_{i+1})}+\mathbf
{1}_{[a_{n},\infty)}.
\end{eqnarray*}
Then $F_{1}$, $F_{2}$ and $\theta$ satisfy all the requirements.
\end{pf}

\begin{coro}\label{nby2}
For any $F^{*}\in\mathcal{D}_{n}$, $n=2,3,\ldots,$ there exist
$F_{k}\in\mathcal{D}_{2}$ and $\theta_{k}\in[0,1]$, $k=1,2,\ldots
,l$, such that
%
\[
F^{*}=\sum_{k=1}^{l}\theta_{k}F_{k},\qquad \sum_{k=1}^{l}\theta_{k}=1.
\]
\end{coro}

\begin{pf}
Since $F^{*}\in\mathcal{D}_{n}\subseteq \mathcal{S}_{n}$, it
follows immediately from Lemma~\ref{2by1} that there exist $F_{k}\in
\mathcal{S}_{2}$ and $\theta_{k}\in[0,1]$, $k=1,2,\ldots,l$, such that
%
\begin{eqnarray*}
F^{*}=\sum_{k=1}^{l}\theta_{k}F_{k},\qquad \sum_{k=1}^{l}\theta_{k}=1,\qquad
\int_{0}^{\infty}\bigl(1-F_{k}(x)\bigr)\,\dd x=\int_{0}^{\infty}\bigl(1-F^{*}(x)\bigr)\,\dd x
\end{eqnarray*}
for $k=1,2,\ldots,l.$
Because $F^{*}\in\mathcal{D}_{n}\subseteq \mathcal{D}$, we have
%
\[
\int_{0}^{\infty}\bigl(1-F_{k}(x)\bigr)\,\dd x=\int_{0}^{\infty}\bigl(1-F^{*}(x)\bigr)\,\dd
x\leqslant s,
\]
which implies that $F_{k}\in \mathcal{D}$. Therefore, $F_{k}\in
\mathcal{S}_{2}\cap\mathcal{D}=\mathcal{D}_{2}$.
\end{pf}

\begin{pf*}{Proof of Lemma~\ref{f=f2}}
Suppose for some $F_{0}\in\mathcal{D}$, we have\break
$
\sup_{F\in\mathcal{D}_{2}}J_{D}(F)<J_{D}(F_{0})<\infty.
$
Construct a sequence of step CDFs, $F_{m}, m=1,2,\ldots,$ satisfying
$F_{m}\geqslant F_0$ and $\lim_{m\rightarrow\infty}F_m(x)=F_0(x)$
a.e. Clearly $F_{m}\in\mathcal{D}$, and it follows from the
dominated convergence theorem that\break $\lim_{m\rightarrow\infty
}J_D(F_{m})=J_D(F_0)$.
So there exists $F^{*}\in\mathcal{D}_{n}$ for some $n\geqslant 2$
such that
$J_{D}(F^{*})>  \sup_{F\in\mathcal{D}_{2}}J_{D}(F). $
By Corollary~\ref{nby2}, there exist $\bar{F}_{k}\in\mathcal
{D}_{2}$ and $\theta_{k}\in[0,1]$, $k=1,2,\ldots,l$, such that
%
\[
F^{*}=\sum_{k=1}^{l}\theta_{k}\bar{F}_{k},\qquad \sum
_{k=1}^{l}\theta_{k}=1.
\]
However, recalling that $w$ is convex, we have
\[
J_{D}(F^{*})=J_{D}\Biggl(\sum_{k=1}^{l}\theta_{k}\bar
{F}_{k}\Biggr)\leqslant\sum_{k=1}^{l}\theta_{k}J_{D}(\bar
{F}_{k})\leqslant\sup_{F\in\mathcal{D}_{2}}J_{D}(F),
\]
which leads to a contradiction.
\end{pf*}
\section{\texorpdfstring{Proof of Lemma \lowercase{\protect\ref{g=g2}}}{Proof of Lemma 4.4}}\label{AD}

Suppose
%
\begin{equation}\label{assu}
\sup_{G\in\mathcal{Q}}J_{Q}(G)>\sup_{G\in\mathcal{Q}_{2}}J_{Q}(G).
\end{equation}
By the monotone convergence theorem, we can find a sequence of
essentially bounded quantile functions $G_n\in\mathcal{Q}$,
$n=1,2,\ldots,$ so that
$\lim_{n\rightarrow\infty}J_Q(G_n)=\sup_{G\in\mathcal
{Q}}J_{Q}(G)$. For each fixed $n$,
by the dominated convergence theorem, there is a sequence of step
functions $G_{n,k}\in\mathcal{Q}$
with $\lim_{k\rightarrow\infty}J_Q(G_{n,k})=J_Q(G_n)$.
So we can find a step function $G_{0}\in\mathcal{Q}$, written as
\begin{eqnarray}
G_{0}(x)=a_{0}+\sum_{i=1}^{n}b_{i}\mathbf{1}_{(c_{i},1]}(x),
\nonumber\\
\eqntext{ a_{0}>
0,  b_{i}> 0,  0<c_{1}<\cdots<c_{n}<1, x\in(0,1),}
\end{eqnarray}
such that
$J_{Q}(G_{0})>\sup_{G\in\mathcal{Q}_{2}}J_{Q}(G).$
Since $G_{0}\in\mathcal{Q}$, we have
\[
\bar{s}:=a_{0}+\sum_{i=1}^{n}b_{i}(1-c_{i})\equiv\int
_{0}^{1}G_0(x)\,\dd x\leqslant s.
\]
Let $0<\varepsilon<a_{0}$. Set
%
\begin{eqnarray*}
a_{i}=\varepsilon,\qquad
\alpha_{i}:=\frac{b_{i}(1-c_{i})}{\bar{s}-\varepsilon}> 0,\qquad
i=1,\ldots,n,\qquad \alpha_{n+1}:= \frac{a_{0}-\varepsilon}{\bar
{s}-\varepsilon}>0,\\
G_{i}(x):=a_{i}+\frac{b_{i}}{\alpha_{i}}\mathbf{1}_{(c_{i},1]}(x),\qquad
i=1,\ldots,n, \qquad G_{n+1}(x):=\bar{s} \qquad\forall x\in(0,1).
\end{eqnarray*}
It is easy to check that $ G_{i}\in\mathcal{Q}_{2}$, $i=1,\ldots
,n+1$, and
\begin{eqnarray*}
G_{0}(x)=\sum_{i=1}^{n+1}\alpha_{i}G_{i}(x),\qquad \sum_{i=1}^{n+1}\alpha
_{i}=1 \qquad \forall x\in(0,1).
\end{eqnarray*}
Recalling that $u$ is convex, we have
\begin{eqnarray*}
\sup_{G\in\mathcal{Q}_{2}}J_{Q}(G)<J_{Q}(G_{0})=J_{Q}\Biggl(\sum
_{i=1}^{n+1}\alpha_{i}G_{i}\Biggr)\leqslant\sum_{i=1}^{n+1}\alpha
_{i}J_{Q}(G_{i})\leqslant\sup_{G\in\mathcal{Q}_{2}}J_{Q}(G),
\end{eqnarray*}
which is a contradiction. So \eqref{assu} is false and the proof is complete.

\section{\texorpdfstring{Proof of Theorem \lowercase{\protect\ref{uconcavews}}}{Proof of Theorem 5.2}}\label{AE}
The key idea of this proof is to show that one needs only to search
among a special class of quantile functions in order to solve the
relaxed problem \eqref{lagr}.
To this end,
fix $G\in\mathcal{Q}$ and $\lambda\geqslant0$, and let
\[
x_{0}:=\sup\biggl\{x\in(0,q] \Big| G(x)\leqslant(u')^{-1}_{l}\biggl(\frac{\lambda
}{w'(1-x)}\biggr) \biggr\}\vee0.
\]

If $x_0>0$, we define
\[
x_{1}=\sup\biggl\{x\in(q,1) \Big| (u')^{-1}_{l}\biggl(\frac{\lambda
}{w'(1-x)}\biggr)\leqslant G(x_{0}) \biggr\}\vee q
\]
and
\[
\hat{G}_{\lambda}(x)=G(x_{0})\mathbf
{1}_{(0,x_{1}]}(x)+(u')^{-1}_{l}\biggl(\frac{\lambda}{w'(1-x)}\biggr)\mathbf
{1}_{(x_{1},1)}(x) \qquad \forall x\in[0,1).
\]
Then $\hat{G}_{\lambda}$ is also a quantile function. We now show
that $J^{\lambda}_{Q}(G)\leqslant J^{\lambda}_{Q}(\hat{G}_{\lambda})$.
Noting $(u')^{-1}_{l}(\frac{\lambda}{w'(1-x)})$ is nonincreasing in
$x\in(0,x_{0})$, we deduce
\begin{eqnarray}
G(x)\leqslant G(x_{0})=G(x_{0}-)\leqslant (u')^{-1}_{l}\biggl(\frac{\lambda
}{w'((1-x_{0})+)}\biggr)\leqslant
(u')^{-1}_{l}\biggl(\frac{\lambda}{w'(1-x)}\biggr)\nonumber\\
\eqntext{\forall x\in(0,x_{0}).}
\end{eqnarray}
Since $f^{\lambda}(x,\cdot)$ is nondecreasing on $[0,
(u')^{-1}_{l}(\frac{\lambda}{w'(1-x)})]$ when $x\in(0,x_{0})$, we have
\begin{eqnarray*}
f^{\lambda}(x,G(x)) \leqslant f^{\lambda}(x,G(x_{0}))=f^{\lambda
}(x,\hat{G}_{\lambda}(x)) \qquad\forall x\in(0,x_{0}).
\end{eqnarray*}
Next, for any $x\in(x_{0},x_{1})$,
$G(x)\geqslant G(x_{0}) \geqslant (u')^{-1}_{l}(\frac{\lambda}{w'(1-x)})$
and
$f^{\lambda}(x,\cdot)$ is nonincreasing on $[ (u')^{-1}_{l}(\frac
{\lambda}{w'(1-x)}),\infty)$.
Hence,
\[
f^{\lambda}(x,G(x)) \leqslant f^{\lambda}(x,G(x_{0}))=f^{\lambda
}(x,\hat{G}_{\lambda}(x))\qquad \forall x\in(x_{0},x_{1}).
\]
Finally,
\begin{eqnarray*}
f^{\lambda}(x,G(x)) \leqslant f^{\lambda}\biggl(x,(u')^{-1}_{u}\biggl(\frac
{\lambda}{w'(1-x)}\biggr)\biggr)=f^{\lambda}(x,\hat{G}_{\lambda}(x))\qquad \forall
x\in(x_{1},1).
\end{eqnarray*}
Therefore,
\[
J^{\lambda}_{Q}(G)=\int_{0}^{1}f^{\lambda}(x,G(x))\,\dd x\leqslant
\int_{0}^{1}f^{\lambda}(x,\hat{G}_{\lambda}(x))\,\dd x=J^{\lambda
}_{Q}(\hat{G}_{\lambda}).
\]

If $x_0=0$, we define
\[
x_{1}=\sup\biggl\{x\in(q,1) \Big| (u')^{-1}_{l}\biggl(\frac{\lambda
}{w'(1-x)}\biggr)\leqslant G(0+) \biggr\}\vee q
\]
and
\begin{eqnarray*}
\hat{G}_{\lambda}(x)&=&G(0+)\mathbf
{1}_{(0,x_{1}]}(x)+(u')^{-1}_{l}\biggl(\frac{\lambda}{w'(1-x)}\biggr)\mathbf
{1}_{(x_{1},1)}(x) \qquad\forall x\in[0,1).
\end{eqnarray*}
%
A similar argument as above shows that $J^{\lambda}_{Q}(G)\leqslant
J^{\lambda}_{Q}(\hat{G}_{\lambda})$.

We have now proved that in order to find an optimal quantile function
one needs only to consider functions of the form
$G(x)=a\mathbf{1}_{(0,q]}(x)+a\vee(u')^{-1}_{l}(\frac{\lambda
}{w'(1-x)})\mathbf{1}_{(q,1)}(x)$,
where the parameters $a$ and $\lambda$ are subject to the constraints
in \eqref{para}.
Note that the last equality constraint in \eqref{para} was due to the
fact that the following payoff function is nondecreasing in $a$.
The payoff under the above $G$ is
\begin{eqnarray*}
J(a,\lambda)&:=&J_{Q}(G)\\
&=&\bigl(1-w(1-q)\bigr)u(a)+\int_{q}^{1} u\biggl( a\vee
(u')_{l}^{-1}\biggl(\frac{\lambda}{w'(1-x)}\biggr)\biggr)w'(1-x)\,\dd x,
\end{eqnarray*}
which is exactly the objective function of \eqref{para}. Since the
optimal solution $a^*>0$, the corresponding $G^*$ defined in \eqref
{sshape} is a quantile. The proof is complete.
\end{appendix}

\section*{Acknowledgments}
We are grateful for comments from seminar and conference participants
at Oxford, ETH, University of Amsterdam, Chinese Academy of Sciences,
University of Hong Kong, Chinese University of Hong Kong, Carnegie
Mellon, University of Alberta, Fudan, McMaster, the 2009 Workshop on
Optimal Stopping
and Singular Stochastic Control Problems in Finance in Singapore, the
1st Columbia--Oxford Joint Workshop in Mathematical Finance in New
York, the 6th World Congress of the Bachelier Finance Society in
Toronto, the 2011 Conference on Modeling and Managing Financial Risks
in Paris, the 2011 Workshop on Recent
Developments in Mathematical Finance in Stockholm, the 2001 Conference
on Stochastic Analysis and Applications in Financial Mathematics in
Beijing and the 2011 International Workshop on Finance in Kyoto. We
thank Jan Ob{\l}{\'o}j for many helpful discussions on the Skorokhod
embedding problem, as well as the two anonymous referees for their
comments that have led to a much improved version of the paper.

%

%

%

%


\printaddresses


\begin{thebibliography}{36}

\bibitem{AY}
\begin{bincollection}[mr]
\bauthor{\bsnm{Az{\'e}ma},~\bfnm{Jacques}\binits{J.}} \AND
  \bauthor{\bsnm{Yor},~\bfnm{Marc}\binits{M.}}
(\byear{1979}).
\btitle{Une solution simple au probl\`eme de {S}korokhod}.
In \bbooktitle{S\'eminaire de {P}robabilit\'es, {XIII} ({U}niv. {S}trasbourg,
  {S}trasbourg, 1977/78)}.
\bseries{Lecture Notes in Math.}
\bvolume{721}
\bpages{90--115}.
\bpublisher{Springer}, \baddress{Berlin}.
\bid{mr={0544782}}
\bptok{imsref}%
\end{bincollection}
\endbibitem

\bibitem{B}
\begin{bmisc}[auto:STB|2012/06/04|06:16:18]
\bauthor{\bsnm{Barberis},~\bfnm{N.}\binits{N.}}
(\byear{2012}).
\bhowpublished{A model of casino gambling. \textit{Management Science}
\textbf{58} 35--51}.
\bptok{imsref}%
\end{bmisc}
\endbibitem

\bibitem{BMZ}
\begin{bmisc}[auto:STB|2012/06/04|06:16:18]
\bauthor{\bsnm{Bj{\"o}rk},~\bfnm{T.}\binits{T.}},
  \bauthor{\bsnm{Murgoci},~\bfnm{A.}\binits{A.}} \AND
  \bauthor{\bsnm{Zhou},~\bfnm{X.~Y.}\binits{X.~Y.}}
  (\byear{2012}).
\bhowpublished{Mean--variance portfolio optimization with state
  dependent risk aversion. \textit{Math. Finance}. To appear}.
\bptok{imsref}%
\end{bmisc}
\endbibitem


\bibitem{CD}
\begin{barticle}[mr]
\bauthor{\bsnm{Carlier},~\bfnm{G.}\binits{G.}} \AND
  \bauthor{\bsnm{Dana},~\bfnm{R.~A.}\binits{R.~A.}}
(\byear{2005}).
\btitle{Rearrangement inequalities in non-convex insurance models}.
\bjournal{J. Math. Econom.}
\bvolume{41}
\bpages{483--503}.
\bid{doi={10.1016/j.jmateco.2004.12.004}, issn={0304-4068}, mr={2143822}}
\bptok{imsref}%
\end{barticle}
\endbibitem

\bibitem{CMM}
\begin{barticle}[mr]
\bauthor{\bsnm{Castagnoli},~\bfnm{Erio}\binits{E.}},
  \bauthor{\bsnm{Maccheroni},~\bfnm{Fabio}\binits{F.}} \AND
  \bauthor{\bsnm{Marinacci},~\bfnm{Massimo}\binits{M.}}
(\byear{2004}).
\btitle{Choquet insurance pricing: A caveat}.
\bjournal{Math. Finance}
\bvolume{14}
\bpages{481--485}.
\bid{doi={10.1111/j.0960-1627.2004.00201.x}, issn={0960-1627}, mr={2070175}}
\bptok{imsref}%
\end{barticle}
\endbibitem

\bibitem{DANA}
\begin{barticle}[mr]
\bauthor{\bsnm{Dana},~\bfnm{Rose-Anne}\binits{R.-A.}}
(\byear{2005}).
\btitle{A representation result for concave {S}chur concave functions}.
\bjournal{Math. Finance}
\bvolume{15}
\bpages{613--634}.
\bid{doi={10.1111/j.1467-9965.2005.00253.x}, issn={0960-1627}, mr={2168522}}
\bptok{imsref}%
\end{barticle}
\endbibitem

\bibitem{DP}
\begin{bbook}[auto:STB|2012/06/04|06:16:18]
\bauthor{\bsnm{Dixit},~\bfnm{A.}\binits{A.}} \AND
  \bauthor{\bsnm{Pindyck},~\bfnm{R.}\binits{R.}}
(\byear{1994}).
\btitle{Investment Under Uncertainty}.
\bpublisher{Princeton Univ. Press}, \baddress{Princeton}.
\bptok{imsref}%
\end{bbook}
\endbibitem

\bibitem{EL}
\begin{bmisc}[auto:STB|2012/06/04|06:16:18]
\bauthor{\bsnm{Ekeland},~\bfnm{I.}\binits{I.}} \AND
\bauthor{\bsnm{Lazrak},~\bfnm{A.}\binits{A.}}
(\byear{2006}).
\bhowpublished{Being serious about non-commitment: Subgame perfect
  equilibrium in continuous time. Working paper}.
\bptok{imsref}%
\end{bmisc}
\endbibitem

\bibitem{Fri}
\begin{bbook}[mr]
\bauthor{\bsnm{Friedman},~\bfnm{Avner}\binits{A.}}
(\byear{1975}).
\btitle{Stochastic Differential Equations and Applications, {V}ols. 1--2}.
\bpublisher{Academic Press}, \baddress{New York}.
\bptok{imsref}%
\end{bbook}
\endbibitem

\bibitem{Hall}
\begin{barticle}[mr]
\bauthor{\bsnm{Hall},~\bfnm{W.~J.}\binits{W.~J.}}
(\byear{1969}).
\btitle{Embedding submartingales in {W}iener processes with drift, with
  applications to sequential analysis}.
\bjournal{J. Appl. Probab.}
\bvolume{6}
\bpages{612--632}.
\bid{issn={0021-9002}, mr={0256459}}
\bptok{imsref}%
\end{barticle}
\endbibitem

\bibitem{HZ}
\begin{bmisc}[auto:STB|2012/06/04|06:16:18]
\bauthor{\bsnm{He},~\bfnm{X.~D.}\binits{X.~D.}} \AND
  \bauthor{\bsnm{Zhou},~\bfnm{X.~Y.}\binits{X.~Y.}}
  (\byear{2009}).
\bhowpublished{Hope, fear, and aspiration. Working paper}.
\bptok{imsref}%
\end{bmisc}
\endbibitem

\bibitem{HZ0}
\begin{barticle}[mr]
\bauthor{\bsnm{He},~\bfnm{Xue~Dong}\binits{X.~D.}} \AND
  \bauthor{\bsnm{Zhou},~\bfnm{Xun~Yu}\binits{X.~Y.}}
(\byear{2011}).
\btitle{Portfolio choice via quantiles}.
\bjournal{Math. Finance}
\bvolume{21}
\bpages{203--231}.
\bid{doi={10.1111/j.1467-9965.2010.00432.x}, issn={0960-1627}, mr={2790902}}
\bptok{imsref}%
\end{barticle}
\endbibitem

\bibitem{Vicky}
\begin{bmisc}[auto:STB|2012/06/04|06:16:18]
\bauthor{\bsnm{Henderson},~\bfnm{V.}\binits{V.}}
(\byear{2012}).
\bhowpublished{Prospect theory, partial liquidation and the disposition
  effect. \textit{Management Science} \textbf{58} 445--460}.
\bptok{imsref}%
\end{bmisc}
\endbibitem


\bibitem{Jin&Zhou}
\begin{barticle}[mr]
\bauthor{\bsnm{Jin},~\bfnm{Hanqing}\binits{H.}} \AND
  \bauthor{\bsnm{Zhou},~\bfnm{Xun~Yu}\binits{X.~Y.}}
(\byear{2008}).
\btitle{Behavioral portfolio selection in continuous time}.
\bjournal{Math. Finance}
\bvolume{18}
\bpages{385--426}.
\bid{doi={10.1111/j.1467-9965.2008.00339.x}, issn={0960-1627}, mr={2427728}}
\bptnote{check related}%
\bptok{imsref}%
\end{barticle}
\endbibitem

\bibitem{JXZ}
\begin{barticle}[auto:STB|2012/06/04|06:16:18]
\bauthor{\bsnm{Jin},~\bfnm{Hanqing}\binits{H.}},
\bauthor{\bsnm{Xu},~\bfnm{Z.}\binits{Z.}}
\AND
\bauthor{\bsnm{Zhou},~\bfnm{Xun~Yu}\binits{X.~Y.}}
(\byear{2008}).
\btitle{A convex stochastic optimization problem arising from portfolio selection}.
\bjournal{Math. Finance}
\bvolume{21}
\bpages{775--793}.
\bptok{imsref}%
\end{barticle}
\endbibitem

\bibitem{Kahneman&Tversky}
\begin{barticle}[auto:STB|2012/06/04|06:16:18]
\bauthor{\bsnm{Kahneman},~\bfnm{D.}\binits{D.}} \AND
  \bauthor{\bsnm{Tversky},~\bfnm{A.}\binits{A.}}
(\byear{1979}).
\btitle{Prospect theory: An analysis of decision under risk}.
\bjournal{Econometrica}
\bvolume{46}
\bpages{171--185}.
\bptok{imsref}%
\end{barticle}
\endbibitem

\bibitem{Lopes}
\begin{barticle}[auto:STB|2012/06/04|06:16:18]
\bauthor{\bsnm{Lopes},~\bfnm{L.~L.}\binits{L.~L.}}
(\byear{1987}).
\btitle{Between hope and fear: The psychology of risk}.
\bjournal{Adv. Experimental Social Psychology}
\bvolume{20}
\bpages{255--295}.
\bptok{imsref}%
\end{barticle}
\endbibitem

\bibitem{NO}
\begin{barticle}[auto:STB|2012/06/04|06:16:18]
\bauthor{\bsnm{Nishimura},~\bfnm{K.}\binits{K.}} \AND
  \bauthor{\bsnm{Ozaki},~\bfnm{H.}\binits{H.}}
(\byear{2007}).
\btitle{Irrevsible investment and Knightian uncertainty}.
\bjournal{J.~Econom. Theory}
\bvolume{136}
\bpages{668--694}.
\bptok{imsref}%
\end{barticle}
\endbibitem

\bibitem{Jan}
\begin{barticle}[mr]
\bauthor{\bsnm{Ob{\l}{\'o}j},~\bfnm{Jan}\binits{J.}}
(\byear{2004}).
\btitle{The {S}korokhod embedding problem and its offspring}.
\bjournal{Probab. Surv.}
\bvolume{1}
\bpages{321--390}.
\bid{doi={10.1214/154957804100000060}, issn={1549-5787}, mr={2068476}}
\bptok{imsref}%
\end{barticle}
\endbibitem

\bibitem{Parato}
\begin{bmisc}[auto:STB|2012/06/04|06:16:18]
\bauthor{\bsnm{Parato},~\bfnm{V.}\binits{V.}}
(\byear{1897}).
\bhowpublished{Cours D'\'Economie Politique.
Lausanne and Paris.}
\bptok{imsref}%
\end{bmisc}
\endbibitem

\bibitem{PrelecD98df}
\begin{barticle}[mr]
\bauthor{\bsnm{Prelec},~\bfnm{Drazen}\binits{D.}}
(\byear{1998}).
\btitle{The probability weighting function}.
\bjournal{Econometrica}
\bvolume{66}
\bpages{497--527}.
\bid{doi={10.2307/2998573}, issn={0012-9682}, mr={1627026}}
\bptok{imsref}%
\end{barticle}
\endbibitem


\bibitem{Riedel}
\begin{barticle}[mr]
\bauthor{\bsnm{Riedel},~\bfnm{Frank}\binits{F.}}
(\byear{2009}).
\btitle{Optimal stopping with multiple priors}.
\bjournal{Econometrica}
\bvolume{77}
\bpages{857--908}.
\bid{doi={10.3982/ECTA7594}, issn={0012-9682}, mr={2531363}}
\bptok{imsref}%
\end{barticle}
\endbibitem

\bibitem{Schied1}
\begin{barticle}[mr]
\bauthor{\bsnm{Schied},~\bfnm{Alexander}\binits{A.}}
(\byear{2004}).
\btitle{On the {N}eyman--{P}earson problem for law-invariant risk measures and
  robust utility functionals}.
\bjournal{Ann. Appl. Probab.}
\bvolume{14}
\bpages{1398--1423}.
\bid{doi={10.1214/105051604000000341}, issn={1050-5164}, mr={2071428}}
\bptok{imsref}%
\end{barticle}
\endbibitem

\bibitem{SXZ}
\begin{barticle}[mr]
\bauthor{\bsnm{Shiryaev},~\bfnm{Albert}\binits{A.}},
  \bauthor{\bsnm{Xu},~\bfnm{Zuoquan}\binits{Z.}} \AND
  \bauthor{\bsnm{Zhou},~\bfnm{Xun~Yu}\binits{X.~Y.}}
(\byear{2008}).
\btitle{Thou shalt buy and hold}.
\bjournal{Quant. Finance}
\bvolume{8}
\bpages{765--776}.
\bid{doi={10.1080/14697680802563732}, issn={1469-7688}, mr={2488734}}
\bptok{imsref}%
\end{barticle}
\endbibitem

\bibitem{Shr}
\begin{bbook}[mr]
\bauthor{\bsnm{Shiryayev},~\bfnm{A.~N.}\binits{A.~N.}}
(\byear{1978}).
\btitle{Optimal Stopping Rules}.
\bpublisher{Springer}, \baddress{New York}.
\bid{mr={0468067}}
\bptok{imsref}%
\end{bbook}
\endbibitem

\bibitem{Skorokhod}
\begin{bbook}[auto:STB|2012/06/04|06:16:18]
\bauthor{\bsnm{Skorokhod},~\bfnm{A.~V.}\binits{A.~V.}}
(\byear{1961}).
\btitle{Issledovaniya po Teorii Sluchainykh Protsessov (Stokhas-Ticheskie
Differentsialnye Uravneniya i Predelnye Teoremy dlya Protsessov Markova)}.
\bpublisher{Izdat. Kiev. Univ.}, \baddress{Kiev, Ukraine}.
\bptok{imsref}%
\end{bbook}
\endbibitem

\bibitem{TverskyAFoxC95pd}
\begin{barticle}[auto:STB|2012/06/04|06:16:18]
\bauthor{\bsnm{Tversky},~\bfnm{A.}\binits{A.}}
\AND
\bauthor{\bsnm{Fox},~\bfnm{C. R.}\binits{C. R.}}
(\byear{1995}).
\btitle{Weighing risk and uncertainty}.
\bjournal{Psychological Rev.}
\bvolume{102}
\bpages{269--283}.
\bptok{imsref}%
\end{barticle}
\endbibitem

\bibitem{T&K1992}
\begin{barticle}[auto:STB|2012/06/04|06:16:18]
\bauthor{\bsnm{Tversky},~\bfnm{A.}\binits{A.}} \AND
  \bauthor{\bsnm{Kahneman},~\bfnm{D.}\binits{D.}}
(\byear{1992}).
\btitle{Advances in prospect theory: Cumulative representation of uncertainty}.
\bjournal{J. Risk Uncertainty}
\bvolume{5}
\bpages{297--323}.
\bptok{imsref}%
\end{barticle}
\endbibitem


\bibitem{Wang}
\begin{barticle}[mr]
\bauthor{\bsnm{Wang},~\bfnm{Shaun}\binits{S.}}
(\byear{1995}).
\btitle{Insurance pricing and increased limits ratemaking by proportional
  hazards transforms}.
\bjournal{Insurance Math. Econom.}
\bvolume{17}
\bpages{43--54}.
\bid{doi={10.1016/0167-6687(95)91054-P}, issn={0167-6687}, mr={1363642}}
\bptok{imsref}%
\end{barticle}
\endbibitem

\bibitem{WY}
\begin{barticle}[mr]
\bauthor{\bsnm{Wang},~\bfnm{Shaun~S.}\binits{S.~S.}} \AND
  \bauthor{\bsnm{Young},~\bfnm{Virginia~R.}\binits{V.~R.}}
(\byear{1998}).
\btitle{Risk-adjusted credibility premiums using distorted probabilities}.
\bjournal{Scand. Actuar. J.}
\bvolume{2}
\bpages{143--165}.
\bid{doi={10.1080/03461238.1998.10413999}, issn={0346-1238}, mr={1659301}}
\bptok{imsref}%
\end{barticle}
\endbibitem

\bibitem{Yaari}
\begin{barticle}[mr]
\bauthor{\bsnm{Yaari},~\bfnm{Menahem~E.}\binits{M.~E.}}
(\byear{1987}).
\btitle{The dual theory of choice under risk}.
\bjournal{Econometrica}
\bvolume{55}
\bpages{95--115}.
\bid{doi={10.2307/1911158}, issn={0012-9682}, mr={0875518}}
\bptok{imsref}%
\end{barticle}
\endbibitem

\end{thebibliography}
\end{document}